\documentclass[11pt,reqno,sumlimits]{amsart}
\input diagrams.tex
\usepackage{amssymb,amscd, epsfig, mathrsfs}
\newtheorem{thm}{Theorem}[section]
\newtheorem{cor}[thm]{Corollary}
\newtheorem{lem}[thm]{Lemma}
\newtheorem{prop}[thm]{Proposition}
\newtheorem{df-pr}[thm]{Definition-Proposition}
\newtheorem{conj}[thm]{Conjecture}

\theoremstyle{definition}
\newtheorem{defn}[thm]{Definition}
\newtheorem{rem}[thm]{Remark}

\newtheorem{notation}[thm]{Notation}

\newcommand{\GSI}{{G^I\rtimes\Sigma_I}}

\newcommand{\Q}{{\mathbb Q}}
\newcommand{\C}{{\mathbb C}}

\newcommand{\Z}{{\mathbb Z}}
\newcommand{\h}{{\mathbb H}}

\newcommand{\Pf}{{\bf Proof: }}
\newcommand{\1}{{\bf 1}}

\newcommand{\g}{{\mathfrak  g}}

\renewcommand{\h}{{\mathfrak  h}}
\newcommand{\w}{{\mathfrak  w}}
\newcommand{\f}{{\mathfrak  f}}
\newcommand{\e}{{\mathfrak  e}}
\newcommand{\fa}{{\mathfrak  a}}
\newcommand{\fc}{{\mathfrak c}}

\newcommand{\fu}{\mathfrak{u}}

\newcommand{\cA}{{\mathcal A}}

\newcommand{\cD}{{\mathcal D}}
\newcommand{\cE}{{\mathcal E}}

\newcommand{\cH}{{\mathscr H}}

\newcommand{\cO}{{\mathcal O}}

\newcommand{\bS}{{\mathbf S}}

\newcommand{\cX}{{\mathcal X}}
\newcommand{\cU}{{\mathcal U}}

\newcommand{\cZ}{{\mathcal Z}}
\newcommand{\sS}{{\mathscr S}}
\newcommand{\sR}{{\mathscr R}}
\newcommand{\sL}{{\mathscr L}}
\newcommand{\sO}{{\mathscr O}}
\newcommand{\sM}{{\mathscr M}}

\newcommand{\bu}{{\mathbf u}}
\newcommand{\bm}{{\mathbf m}}

\newcommand{\br}{{\mathbf r}}
\newcommand{\bx}{{\mathbf x}}
\newcommand{\bq}{{\mathbf q}}
\newcommand{\bp}{{\mathbf p}}
\newcommand{\by}{{\mathbf y}}
\newcommand{\bs}{{\mathbf s}}
\newcommand{\bv}{{\mathbf v}}

\newcommand{\bc}{{\mathbf c}}

\newcommand{\be}{{\mathbf e}}

\newcommand{\bgot}{\bigotimes}

\newcommand{\lan}{{\langle}}
\newcommand{\ran}{{\rangle}}

\newcommand{\iso}{\stackrel{\simeq}{\longrightarrow}}

\newcommand{\inc}{\hookrightarrow}

\newcommand{\Supp}{\operatorname{Supp}}

\newcommand{\Hom}{\operatorname{Hom}}

\newcommand{\Map}{\operatorname{Map}}

\newcommand{\Tr}{\operatorname{Tr}}

\newcommand{\SL}{{\mbox{SL}}}

\newcommand{\fbar }{\overline{f}}

\newsavebox{\savepar}

\numberwithin{equation}{section}

\newcounter{labelflag} 
\newcommand{\labelon}{\setcounter{labelflag}{1}}
\newcommand{\Label}[1]{\ifnum\thelabelflag=1\ifmmode
\makebox[0in][l]{\qquad\fbox{\rm#1}} \else
\marginpar{\vspace{0.7\baselineskip} \hspace{-1.1\textwidth}
\fbox{\rm#1}} \fi \fi \label{#1} } \labelon
\begin{document}
\title{Orbifold Cohomology of a Wreath Product Orbifold}
\author{Tomoo Matsumura}
\address{Department of Mathematics and Statistics\\
  Boston University}
\email{mushmt@math.bu.edu {\bf \today}}%
\begin{abstract}   Let $[X/G]$ be an orbifold which is a global quotient
of a compact almost complex manifold $X$ by a finite group $G$.
Let $\Sigma_n$ be the symmetric group on $n$ letters. Their
semidirect product $G^n \rtimes \Sigma_n$ is called the {\em
wreath product} of $G$ and it naturally acts on the $n$-fold product
$X^n$, yielding the orbifold $[X^n/(G^n\rtimes \Sigma_n)]$. Let
$\cH(X^n, G^n\rtimes \Sigma_n)$ be the stringy cohomology \cite{FG,
JKK1} of the $(G^n\rtimes \Sigma_n)$-space $X^n$. When $G$ is
Abelian, we show that the $G^n$-coinvariants of $\cH(X^n, G^n\rtimes
\Sigma_n)$ is isomorphic to the algebra $\cA\{\Sigma_n\}$ introduced
by Lehn and Sorger \cite{L-S}, where $\cA$ is the orbifold cohomology of
$[X/G]$. We also prove that, if $X$ is a projective surface with
trivial canonical class and $Y$ is a crepant resolution of $X/G$,
then the Hilbert scheme of $n$ points on $Y$, denoted by $Y^{[n]}$,
is a crepant resolution of $X^n/(G^n\rtimes \Sigma_n)$.
Furthermore, if $H^{\ast}(Y)$ is isomorphic to
$H^{\ast}_{orb}([X/G])$, then $H^{*}(Y^{[n]})$ is isomorphic to
$H^{\ast}_{orb}([X^n/(G^n \rtimes \Sigma_n)])$. Thus we verify a
special case of the cohomological hyper-K\"{a}hler resolution
conjecture due to Ruan \cite{R}.
\end{abstract}
\maketitle

\section{{\bf Introduction}}
The stringy cohomology $\cH(X,G)$ of an almost complex manifold $X$
with an action of a finite group $G$ was introduced by Fantechi and
G\"{o}ttsche \cite{FG}.  It is a $G$-Frobenius algebra
\cite{FG,JKK1} which is a $G$-equivariant generalization of a
Frobenius algebra. The space of $G$-coinvariants of $\cH(X,G)$ is
isomorphic as a Frobenius algebra to the Chen-Ruan orbifold
cohomology $H_{orb}^{*}([X/G])$ of the orbifold $[X/G]$.

In this section, assume that the coefficient ring for
cohomology is $\C$. Let
$\mathcal{W}$ be an orbifold and $\pi: Y \rightarrow W$ be a
hyper-K\"{a}hler resolution of the coarse moduli space $W$ of
$\mathcal{W}$. Ruan's cohomological hyper-K\"{a}hler resolution
conjecture \cite{R} predicts that the ordinary cohomology ring of
$Y$ is isomorphic to the orbifold cohomology ring of $\mathcal{W}$.
This is a special case of the cohomological crepant resolution
conjecture and the crepant resolution conjecture \cite{R}. These
conjectures have been verified in many cases, cf. \cite{P, BGP}.

Among the examples which support the cohomological hyper-K\"{a}hler
resolution conjecture, the symmetric product is perhaps the most
fascinating. Let $Y$ be a projective surface with trivial canonical
class. The symmetric group on $n$-letters, $\Sigma_n$, naturally
acts on the $n$-fold product $Y^n$ of $Y$. The Hilbert scheme of $n$
points on $Y$, denoted by $Y^{[n]}$, is a hyper-K\"{a}hler
resolution of the quotient space $Y^n/\Sigma_n$ \cite{Be}. Fantechi
and G\"{o}ttsche \cite{FG} showed that the ring of
$\Sigma_n$-coinvariants of $\cH(Y^n, \Sigma_n)$ is isomorphic to
$H^*(Y^{[n]})$. Their proof proceeds by showing that $\cH(Y^n,
\Sigma_n)$ is isomorphic to the algebra $\cA\{S_n\}$ defined by Lehn
and Sorger \cite{L-S} where $\cA$ is the ordinary cohomology of $X$,
$i.e.$
\[
\cH(Y^n, \Sigma_n) \cong H^*(Y)\{\Sigma_n\}  \ \ \Longrightarrow \ \
H^*_{orb}([Y^n/\Sigma_n]) \cong H^*(Y)\{\Sigma_n\}^{\Sigma_n} \cong
H^*(Y^{[n]})
\]
where the last isomorphism is due to \cite{L-S} (see also
\cite{Ur,QW,LQW}).

In this paper, we consider a generalization of the algebra isomorphism on
the left-hand side of the arrow above. The symmetric group
$\Sigma_n$ naturally acts on the $n$-fold product $G^n$ and their
semidirect product $G^n \rtimes \Sigma_n$ is called the {\em wreath
product} of $G$. It naturally acts on the $n$-fold product $X^n$,
yielding the orbifold $[X^n/(G^n\rtimes \Sigma_n)]$. This orbifold is
called the {\em wreath product orbifold} of a $G$-space $X$.
The linear structure of the orbifold cohomology of a wreath product orbifold
has been studied in a sequence of papers by Qin, Wang and Zhou,
cf. \cite{QW,W,WZ} through a careful analysis of the fix point loci.
However, one of the goals of this paper is to analyze the multiplication
in stringy cohomology and in Chen-Ruan orbifold cohomology of a
wreath product orbifold. The multiplication
in the special case when $X=\C^2$ and $G$ is a finite subgroup
of $\SL_2(\C)$ has been studied in \cite{EG,QW2}.

The main result of this paper is Theorem \ref{ringiso} which proves
that, when $X$ is compact and $G$ is Abelian, the $G^n$-coinvariants of $\cH(X^n,
G^n\rtimes \Sigma_n)$ is isomorphic as a $\Sigma_n$-Frobenius algebra
to the algebra $\cA\{\Sigma_n\}$ where $\cA$ is the orbifold
cohomology of $[X/G]$, $i.e.$
\[
\cH(X^n, G^n\rtimes \Sigma_n)^{G^n} \cong
H^*_{orb}([X/G])\{\Sigma_n\}.
\]
When $G$ is a trivial group, this isomorphism reduces to the
isomorphism defined by Fantechi and G\"{o}ttsche \cite{FG}.

A key role in this paper is played by the formula (\ref{obstruction
bundle}) for the obstruction bundle of the stringy cohomology, which
is proved by Jarvis, Kaufmann, and Kimura \cite{JKK2}. Since their
definition avoids any construction of complex curves, admissible
covers, or moduli spaces, it greatly simplifies the analysis of the
obstruction bundle and allows us to write the obstruction bundle of
$[X^n/(G^n\rtimes\Sigma_n)]$ in terms of the ones of $[X/G]$.
See Proposition \ref{R(gs,ht)}.

To relate our result to Ruan's conjecture, we need to work in the
algebraic category. We observe (cf. \cite{W}) that, if $X/G$ is an even dimensional
Gorenstein variety and $Y$ is a crepant resolution of $X/G$, then
$Y^n/\Sigma_n$ is a crepant resolution of $X^n/(G^n\rtimes \Sigma_n)$.
Hence, if $Y$ is a projective surface with the trivial canonical
class, then $Y^{[n]}$ is a crepant resolution of
$X^n/(G^n\rtimes\Sigma_n)$, $i.e.$ the composition
\[
Y^{[n]} \longrightarrow Y^n/\Sigma_n \longrightarrow
X^n/(G^n\rtimes\Sigma_n)
\]
is a crepant resolution. Together with Theorem \ref{ringiso}, if
$H^{\ast}(Y) \cong H^{\ast}_{orb}([X/G])$, then we obtain
a verification, in a special case, of the cohomological
hyper-K\"{a}hler resolution conjecture:
\[
H^*_{orb}([X^n/(G^n\rtimes\Sigma_n)]) \cong
H^*_{orb}([X/G])\{\Sigma_n\}^{\Sigma_n} \cong
H^*(Y)\{\Sigma_n\}^{\Sigma_n} \cong H^*(Y^{[n]}).
\]
This conjecture in the special case has been verified
in the case when $X=\C^2$ and
$G$ is a finite subgroup of $\SL_2(\C)$ \cite{EG}.

Our main result, Theorem \ref{ringiso}, fits into a larger framework as
follows. We show (Theorem \ref{thm:L-Frob}) that,
if $\cH$ is a $(K\rtimes L)$-Frobenius
algebra for any semidirect product of
finite groups $K$ and $L$, then the space of $K$-coinvariants of $\cH$
is an $L$-Frobenius algebra. Thus, we
may interpret the $\Sigma_n$-Frobenius algebra obtained by
taking the $G^n$-coinvariants of $\cH(X^n, G^n\rtimes \Sigma_n)$ as
the stringy cohomology of the orbifold $[X/G]^n$ with the action of
$\Sigma_n$. Furthermore, in the case of a global quotient, the
stringy cohomology is obtained as the degree zero $G$-equivariant
Gromov-Witten invariants introduced in \cite{JKK1} for an almost
complex manifold $X$ with an action of a finite group $G$, while the
orbifold cohomology is obtained by the degree zero Gromov-Witten
invariants of the orbifold \cite{CR2, AGV}. The fact that the
$G$-coinvariants of the stringy cohomology of $G$-space $X$ is the
orbifold cohomology of $[X/G]$, also follows from the fact that
orbifold Gromov-Witten invariants of $[X/G]$ is obtained from
$G$-equivariant Gromov-Witten invariants of $G$-space $X$ by taking
its ``$G$-invariants" in the sense of \cite{JKK1}. In particular, if
$G$ is a semidirect product of $K$ and $L$ where $L$ acts on $K$,
there should exist a kind of \emph{$L$-equivariant Gromov-Witten
invariants of an orbifold $[X/K]$ with the action of $L$}, which is
equivalent to the ``$K$-invariants" of the $(K\rtimes L)$-equivariant
Gromov-Witten invariants of the $(K\rtimes L)$-space $X$. That is,
the following diagram should hold.
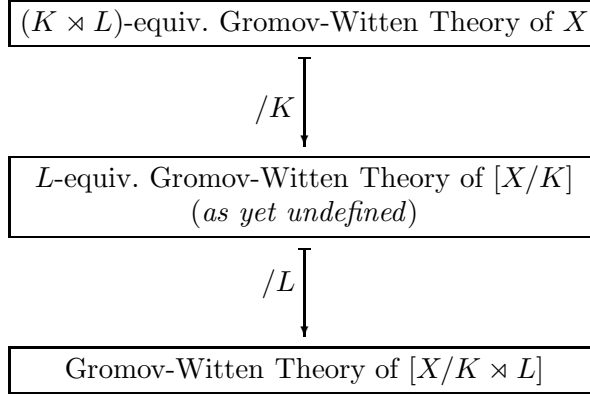
\[
\begin{diagram}
\fbox{\parbox{.50\linewidth}{\centerline{$(K\rtimes L)$-equiv.~Gromov-Witten Theory of $X$}}}
\\
\dMapsto^{/K}
\\
\fbox{\parbox{.50\linewidth}{\centerline{$L$-equiv.~Gromov-Witten Theory of $[X/K]$}
\centerline{(\emph{as yet undefined})}}}
\\
\dMapsto^{/L}
\\
\fbox{\parbox{.50\linewidth}{\centerline{Gromov-Witten Theory of $[X/K\rtimes L]$}}}
\\
\end{diagram}
\]
\

\

%
The structure of the rest of the paper is as follows. In Section \ref{sec:gen},
we review the definition of a $G$-Frobenius algebra and show that,
if $\cH$ is an $(K\rtimes L)$-Frobenius algebra, then
the space of $K$-coinvariants of $\cH$ is an $L$-Frobenius algebra.
In Section \ref{sec:wreath}, we study the
wreath product associated to a finite group $G$. In Sections \ref{LS}
and \ref{Orbifold}, we review the definition of the Lehn-Sorger
algebras $\cA\{\Sigma_n\}$ and prove a geometric formula (Equation
(\ref{61})) for the multiplication in the Lehn-Sorger algebra
associated to $H^*_{orb}([X/G])$. In Section
\ref{sec:wreathproductorbifold}, we prove that there is a canonical
$\Sigma_n$-graded $\Sigma_n$-module isomorphism between
$H^*_{orb}([X/G])\{\Sigma_n\}$ and the space of
$\Sigma_n$-coinvariants of $\cH(X^n, G^n\rtimes \Sigma_n)$. In
Sections \ref{sec:obstruction bundle} and \ref{sec:ring isomorphism},
we compute the obstruction bundle of the stringy cohomology
$\cH(X^n,G^n\rtimes \Sigma_n)$ by using the formula
(\ref{obstruction bundle}) from \cite{JKK2}
and prove, in the case when $G$ is an
Abelian group, that the isomorphism introduced in Section
\ref{sec:wreathproductorbifold} preserves the ring structures. In
Section \ref{sec:ex1}, we work out an example.  In section
\ref{sec:ex2}, we study an example of the simplest case when $G$ is
not Abelian. In Section \ref{sec:crepant}, we verify a special case
of the Ruan's conjecture.

\subsection*{Acknowledgements}
The author is greatly indebted to his thesis advisor Takashi Kimura,
who has provided constant guidance throughout
the course of this project. The author would like to thank Dan Abramovich,
Alastair Craw, Barbara Fantechi, So Okada, Fabio Perroni, Weiqiang Wang
for important advice and useful conversations.

\section{{\bf G-Frobenius algebras and semidirect products}}\label{sec:gen}
Unless otherwise specified, we assume throughout the paper that all
groups are finite and all group actions are right actions. Also,
unless otherwise specified, all of the vector spaces are finite
dimensional and over $\Q$, and all coefficient rings for cohomology
and K-theory are over $\Q$.

We recall the definition of a $G$-Frobenius algebra \cite{JKK1} for
a group $G$.

\begin{defn}
A $G$-graded vector space $\cH:=\bigoplus_{m} \cH_m$ which is
endowed with the structure of a right $G$-module by isomorphisms
$\rho(\gamma) : \cH \iso \cH$ for all $\gamma$ in $G$, is said to be
a \emph{$G$-graded $G$-module} if $\rho(\gamma)$ takes $\cH_m$ to
$\cH_{\gamma^{-1}m\gamma}$ for all $m$ in $G$. We denote a vector in
$\cH_m$ by $v_m$  for any $m \in G$,
\end{defn}

\begin{defn}\label{app}
A tuple $((\cH,\rho),\cdot,\1,\eta)$ is said to be a
\emph{$G$-(equivariant) Frobenius algebra} provided that the
following properties hold:
\begin{itemize}
\item[i)] ($G$-graded $G$-module) $(\cH , \rho)$ is
a $G$-graded $G$-module.
\item[ii)] (Self-invariance) For all $\gamma$ in $G$,
$\rho(\gamma):\cH_{\gamma} \rightarrow \cH_{\gamma}$ is the identity
map.
\item[iii)] (Metric) $\eta$ is a symmetric, non-degenerate,
bilinear form on $\cH$ such that $\eta(v_{m_1} , v_{m_2} )$ is
nonzero only if $ m_1m_2 = 1$.
\item[iv)] ($G$-graded Multiplication) The binary product
$(v_1,v_2) \stackrel{\mu}{\mapsto} v_1\cdot v_2$, called the
\emph{multiplication} on $\cH$, preserves the $G$-grading (
$i.e.$ the multiplication takes $\cH_{m_1} \otimes \cH_{m_2}$ to
$\cH_{m_1m_2}$ ) and is distributive over addition.
\item[v)] (Associativity) The multiplication is associative, $i.e.$
\[
(v_1\cdot v_2) \cdot v_3=v_1\cdot( v_2 \cdot v_3)
\]
 for all $ v_1, v_2,$ and $v_3$ in $\cH$.
\item[vi)] (Braided Commutativity)
The multiplication is invariant with respect to the braiding,
\[
v_{m_1} \cdot v_{m_2} =\left(\rho(m_1^{-1}) v_{m_2}\right) \cdot
v_{m_1}
\]
 for all $m_i \in G$ and all $v_{m_i} \in \cH_{m_i}$ with $i = 1, 2$.
\item[vii)] ($G$-equivariance of the Multiplication)
\[
\rho(\gamma)v_{1}\cdot\rho(\gamma)v_{2}=\rho(\gamma)(v_{1} \cdot
v_{2} )
\]
 for all $\gamma$ in $G$, and all $v_1, v_2 \in \cH$.
\item[viii)] ($G$-invariance of the Metric)
\[
\eta(\rho(\gamma)v_{1} ,\rho(\gamma)v_{2} )=\eta(v_{1},v_{2})
\]
 for all $\gamma$ in $G$, and all $v_1, v_2 \in \cH$.
\item[ix)] (Invariance of the Metric)
\[
\eta(v_1\cdot v_2, v_3)=\eta(v_1, v_2 \cdot v_3)
\]
 for all $v_1, v_2, v_3 \in \cH$.
\item[x)] ($G$-invariant Identity) The element $\1$ in $\cH_1$
is the identity element of the multiplication, which satisfies
$\rho(\gamma)\1 = \1 $ for all $\gamma$ in $G$.
\item[xi)] (Trace Axiom) For all $a, b$ in $G$ and $v$ in
$\cH_{[a,b]}$, if $L_v$ denotes the left multiplication by $v$, then the
following equation is satisfied:
\[
\Tr_{\cH_a} (L_v\rho(b^{-1})) = \Tr_{ \cH_b} (\rho(a)L_v).
\]
\end{itemize}
\end{defn}
\begin{defn}\label{Q-grading}
A $G$-Frobenius algebra $\cH$ is said to be \emph{$\Q$-graded} if we
can write
\[
\cH=\bigoplus_{r \in \Q}\cH_r
\]
and there exists a non-negative integer $d$ such that $\cH_r=0$ if
$r < 0$ or $r>2d$. Furthermore, the $G$-action, $G$-grading,
multiplication respect the $\Q$-grading and the metric has grading
$-2d$. In this paper, we assume that all $G$-Frobenius algebras are
$\Q$-graded.
\end{defn}

\begin{defn}
A $G$-Frobenius algebra when $G=\{1\}$ is called a \emph{Frobenius
algebra}.
\end{defn}

\begin{rem}
We can also define a $G$-Frobenius superalgebra \cite{Ka} by
introducing $\Z/2\Z$-grading and by introducing signs in the
usual manner.
\end{rem}

Let $K$ and $L$ be groups. Let $L$ act on $K$ and we denote the
action of $l\in L$ on $k\in K$ by $ k \stackrel{l}{\mapsto} k^{l}$.
Let $K \rtimes L$ be a semidirect of groups $K$ and $L$ with respect
to this action. We identify $K$ with the normal subgroup $K \rtimes
1$ and hence the adjoint action of $L$ on $K$ can be identified with
the given action of $L$ on $K$, namely, $k^l=l^{-1}kl$.

Let $((\cH,\rho), \cdot,\1 , \eta)$ be a $(K \rtimes L)$-Frobenius
algebra.  Let $\cH_{[a]}:=\oplus_{k \in K} \cH_{ka}$ and let $\pi_K
: \cH \rightarrow \cH$ be the averaging map over $K$:
\[
\pi_K(v) :=\frac{1}{|K|} \sum_{k \in K} \rho(k)v.
\]
The image $\pi_K(\cH)$ is the space of $K$-coinvariants of $\cH$, which we
denote by $\cH^K$. Let $\cH^K_{[l]}:=\pi_K(\cH_{[l]})$.
\begin{thm}\label{thm:L-Frob}
If $\cH$ is a $(K\rtimes L)$-Frobenius algebra, then $\cH^K$ is an
$L$-Frobenius algebra.
\end{thm}
\Pf
All of the properties except the self-invariance property and the
trace axiom follow immediately from those properties of $\cH$. The
self-invariance property of $\cH^K$ is that, for all $l \in L$,
$\rho(l) : \cH^K_{[l]} \rightarrow \cH^K_{[l]}$ is the identity map.
This is true because of the self-invariance property of $\cH$.
Indeed, for all $kl \in K\rtimes L$, $\rho(kl)$ restricted to
$\cH_{kl}$ is the identity map so that $\rho(k)=\rho(l^{-1})$ on
$\cH_{kl}$. Hence, for all $v \in \cH_{k_0l}$,
\[
\rho(l)\pi_K(v)=\frac{1}{|K|}\sum_{k'\in K}
\rho(k')\rho(l)v=\frac{1}{|K|}\sum_{k'\in K}
\rho(k')\rho(k_0^{-1})v=\frac{1}{|K|}\sum_{k''\in
K}\rho(k'')v=\pi_K(v)\] where we have used the self-invariance
property of $\cH$ at the second equality and the third equality is
obtained by the change of variables $k''=k_0k'$. Since any element
of $\cH_{[l]}^K$ is represented by $\pi_{K}(v)$ for some $v \in
\cH_{k_0l}$, we have proved the self-invariance property of $\cH^K$.

The trace axiom for $\cH^{K}$ is the following equality,
\[
\Tr_{\cH^K_{[l_1]}} (L_{v_m} \circ
\rho(l_2^{-1}))=\Tr_{\cH^K_{[l_2]}} (\rho(l_1) \circ L_{v_m}),
\]
for $l_1,l_2 \in L$ and $v_m \in \cH^K_{[m]}$ where $m=[l_1,l_2]$.
 The left-hand side is
\begin{eqnarray*}
\Tr_{\cH^K_{[l_1]}} (L_{v_m} \circ \rho(l_2^{-1}))
&=&\Tr_{\cH_{[l_1]}}(L_{v_m} \circ \rho(l_2^{-1}) \circ \pi_K) \\
&=&\frac{1}{|K|}\sum_{k_1,k} \Tr_{\cH_{k_1l_1}}(L_{v_m} \circ
\rho(l_2^{-1})\circ \rho(k))\\
&=&\frac{1}{|K|}\sum_{k_1,k_2} \Tr_{\cH_{k_1l_1}}(L_{v_m} \circ
\rho(l_2^{-1}k_2^{-1})\\
&=&\frac{1}{|K|}\sum_{k_1,k_2} \Tr_{\cH_{k_2l_2}}(\rho(k_1 l_1)
\circ L_{v_m}),
\end{eqnarray*}
where the third equality is obtained by replacing the parameter
$k^{l_2^{-1}}$ by $k_2^{-1}$ and the fourth equality follows from
the trace axiom for $\cH$. The right-hand side is
\begin{eqnarray*}
\Tr_{\cH^K_{[l_2]}} (\rho(l_1) \circ L_{v_m})
&=&\frac{1}{|K|}\sum_{k_1,k_2}\Tr_{\cH_{k_2l_2}}(\rho(l_1) \circ
L_{v_m}
\circ \rho(k_1)) \\
&=&\frac{1}{|K|}\sum_{k_1,k_2}\Tr_{\cH_{k_2l_2}}(\rho(k_1l_1) \circ
L_{v_m}),
\end{eqnarray*}
where the second equality follows from the cyclicity of the trace
and by replacing the parameter $k_1^{l_1^{-1}}$ by $k_1$. Thus, the
trace axiom holds for the $K$-coinvariants $\cH^K$. \qed

\section{{\bf The wreath product $\GSI$}}\label{sec:wreath}
In this section, we review the wreath product of a group $G$ (cf. \cite{W})
to fix the notation and also to establish a technical lemma
which we will use later.

\begin{notation}
The set of conjugacy classes of $G$ is denoted
by $\overline{G}$. For all $\alpha\in G$, let $\cZ_G(\alpha)$ be the
centralizer of $\alpha$ in $G$. The subgroup generated by the subset
$\{\alpha_k\}_{k=1,\cdots,r}$ of $G$ is denoted by
$\lan\alpha_1,\cdots,\alpha_r\ran$. For a finite set $J$, let $G^J$
be the set of maps, $\Map(J,G)$, from $J$ to $G$. It is, of course,
non-canonically isomorphic as a set to the $|J|$-fold product
$G^{|J|}$ where $|J|$ is the cardinality of $J$. For all $g \in
G^J$, $g_i$ denotes the image of $i\in J$ under $g$ and is called
the \emph{$i$-th component} of $g$. Let $\Delta^J : G \rightarrow
G^J$ be the diagonal map and let $\Delta^J_G$ be the image of $G$
under $\Delta^J$. The same notation is applied to any set, $i.e.$
if $X$ is a manifold, then $X^J:=\Map(J,X)$ and $x_i$ is the
$i$-th component of $x \in X^J$ for all $i \in J$. $\Delta^J_{X}$ is
the image of $X$ under the diagonal map $\Delta^J:X \rightarrow
X^J$.

Let $I$ be a finite set of cardinality $n$ and let $\Sigma_I$ be the
permutation group of the set $I$. For all $\sigma, \tau \in
\Sigma_I$, let $o(\sigma)$ be the set of orbits in $I$ under the
action of the subgroup $\lan\sigma\ran$ and let  $o(\sigma,\tau)$ be
the set of orbits in $I$ under the action of the subgroup
$\lan\sigma,\tau\ran$. Using the natural action of $\Sigma_I$ on
$G^I$, we obtain the semidirect product $G^I \rtimes \Sigma_I$.
Namely, for all $\sigma \in \Sigma_I$, define $g^{\sigma}$ in $G^I$
by $g^{\sigma}_i:= g_{\sigma(i)}$. We denote an element of
$G^I\rtimes \Sigma_I$ by $g\sigma$ for all $g \in G^I$ and $\sigma
\in \Sigma_I$. The product of $g\sigma$ and $h\tau$ in $\GSI$ is
\[
g\sigma\cdot h\tau = gh^{\sigma^{-1}}\sigma \tau
\]
for all $g,h\in G^I$ and $\sigma, \tau \in \Sigma_I$. We also
observe that the action of $G^I$ by conjugation on $\GSI$ preserves
the coset $G^I \sigma=\{g\sigma \ |\ g \in G^I\}$ for each $\sigma
\in \Sigma_I$.
\end{notation}


\begin{defn}For each $a \in o(\sigma)$, choose a representative $i_a \in
a$. For all $g \in G^I$ and $a \in o(\sigma)$, define
\begin{equation}\label{CycleProduct}
\psi^{\sigma}(g)_a:=\prod_{k=0}^{|a|-1}
\g_{\sigma^{|a|-1-k}(i_a)}:=g_{\sigma^{|a|-1}(i_a)}
g_{\sigma^{|a|-2}(i_a)} \cdots g_{\sigma^0(i_a)}
\end{equation}
and let $\psi^{\sigma}(g)$ be the element of $G^{o(\sigma)}$ that
has components $\{\psi^{\sigma}(g)_a\}_{a \in o(\sigma)}$. Call
$\psi^{\sigma}(g)$ a \emph{cycle product} of $g$ with respect to
$\sigma$ associated to $\{i_a\}$.
\end{defn}

For example, let $I=\{1,2,3,4,5\}$ and define $\sigma$ by
$\sigma(1)=3,\sigma(2)=1, \sigma(3)=2, \sigma(4)=5,\sigma(5)=4$. We
have $o(\sigma)=\{a_1, a_2\}$ where $a_1:=\{1,2,3\}$ and
$a_2:=\{4,5\}$. If we choose $2 \in a_1$ and $4 \in a_2$, then
$\psi^{\sigma}(g) \in G^{\{a_1,a_2\}}$ is
\[
\psi^{\sigma}(g)_{a_1}=g_3g_1g_2, \ \ \
\psi^{\sigma}(g)_{a_2}=g_5g_4.
\]
The cycle product depends on the choice of representatives $\{i_a\}$
and if we choose different representatives, then each
$\psi^{\sigma}(g)_a$ will be conjugated by some element $\gamma_a$
in $G$. In the example, if we choose $3 \in a_1$ instead of $2$ and
$5 \in a_2$ instead of $4$, then $\psi^{\sigma}(g)$ is
\[
\psi^{\sigma}(g)_{a_1}=g_1g_2g_3, \ \ \
\psi^{\sigma}(g)_{a_2}=g_4g_5,
\]
and $\gamma_{a_{1}}=g_3$ and $\gamma_{a_2}=g_5$. Hence, the
componentwise $G$-conjugacy class, $\overline{\psi^{\sigma}(g)} \in
\overline{G}^{o(\sigma)}$, is independent of the choice of
representatives $\{i_a\}$.

Now we compute the orbits of the action of $G^I$ by conjugation on
$G^I\rtimes\Sigma_I$.

\begin{defn}
For all $\g \in G^{o(\sigma)}$, let $\overline{\g}$ be the element
in $\overline{G}^{o(\sigma)}$ such that $\overline{\g}_a$ is the
$G$-conjugacy class containing $\g_a$ for each $a \in o(\sigma)$.
Define
\begin{equation}\label{orbit}
\cO_{\g} \sigma :=\left\{\left. g\sigma \in G^I \sigma\ \right|\
\overline{\ \psi^{\sigma}(g)}=\overline{\g}\ \right\}.
\end{equation}
$\cO_{\g} \sigma$ is independent of the choice of representatives
$\{i_a\}$ and clearly, $\cO_{\g}\sigma=\cO_{\h}\sigma$ if and only
if $\overline{\g}=\overline{\h}$.
\end{defn}

\begin{prop}\label{prop:orb}
If $\g=\psi^{\sigma}(g)$, then $\cO_{\g}\sigma$ is the orbit of $g\sigma$
under the action of $G^I$ by conjugation on $G^{I}\rtimes\Sigma_I$.
\end{prop}
\Pf Choose a representative $i_a$ in $a$ for each $a \in o(\sigma)$.
We define two elements $\epsilon_{\g}$ and $\nu^{\sigma}(g)$ in
$G^I$. For each $\g \in G^{o(\sigma)}$ and $j\in I$, let
\begin{equation}\label{def-epsilon}
(\epsilon_{\g})_j:=
\begin{cases}
\g_a & \text{ if } j=i_a \\
1 & \text{ otherwise. }
\end{cases}
\end{equation}
For each $g \in G^I$, let
\begin{equation}\label{nu}
\nu^{\sigma}(g)_{\sigma^m(i_a)}:=g_{\sigma^m(i_a)}
g_{\sigma^{m-1}(i_a)}\cdots g_{\sigma^0(i_a)}
\end{equation}
where $m=0, \cdots, |a|-1$. In particular,
$\nu^{\sigma}(g)_{\sigma^{|a|-1}(i_a)}=\psi^{\sigma}(g)_a$. Note
that each element in $I$ can be represented uniquely by
$\sigma^{m}(i_a)$. If $\g=\psi^{\sigma}(g)$, then we have
\begin{equation}\label{ngn=epsilon}
\nu^{\sigma}(g)^{-1}\cdot g
\sigma\cdot\nu^{\sigma}(g)=\nu^{\sigma}(g)^{-1}\cdot g \cdot
\nu^{\sigma}(g)^{\sigma^{-1}}\sigma=\epsilon_{\g}\sigma
\end{equation}
so that $g\sigma$ and $\epsilon_{\g}\sigma$ are in the same orbit.
 Indeed, if $m\not=0$,
\[
(\nu^{\sigma}(g)^{-1}\cdot g\cdot
\nu^{\sigma}(g)^{\sigma^{-1}})_{\sigma^m(i_a)}
=\nu^{\sigma}(g)_{\sigma^m(i_a)}^{-1}\cdot g_{\sigma^m(i_a)}
\cdot\nu^{\sigma}(g)_{\sigma^{m-1}(i_a)}=1.
\]
If $m=0$, \[ (\nu^{\sigma}(g)^{-1}\cdot
g\cdot\nu^{\sigma}(g)^{\sigma^{-1}})_{i_a}
=\nu^{\sigma}(g)_{\sigma^{|a|-1}(i_a)}=\psi^{\sigma}(g)_a,
\]
since $\sigma^{-1}(i_a)=\sigma^{|a|-1}(i_a)$.

On the other hand, for all $\g$ and $\g'$ in $G^{o(\sigma)}$, there
exists an $f \in G^I$ satisfying $\epsilon_{\g}\sigma
=f^{-1}\epsilon_{\g'}\sigma f$ if and only if there exists $f \in
\prod_{a \in o(\sigma)}\Delta^a_G$ such that $\g_a=f_{i_a}^{-1}\g'_a
f_{i_a}$ for each $a \in o(\sigma)$.

Thus, $g\sigma$ and $g'\sigma$ are in the same orbit if and only if
$\overline{\g}=\overline{\g'}$ where $\g:=\psi^{\sigma}(g)$ and
$\g':=\psi^{\sigma}(g')$. \qed

\begin{rem}
For all $\g \in G^{o(\sigma)}$, we have
\begin{equation}\label{cent}
\cZ_{G^I}(\epsilon_{\g}\sigma)=\prod_{a \in o(\sigma)}
\Delta^a_{\cZ_G(\g_a)}.
\end{equation}
In particular, $\cZ_{G^I}(\epsilon_{\g}\sigma)=\prod_{a \in
o(\sigma)} \Delta^a_G$ if $G$ is Abelian.
\end{rem}
\begin{lem}\label{surj}
Suppose that $G$ is an Abelian group and that $\lan\sigma,\tau\ran$
acts transitively on $I$. Let $\g\in G^{o(\sigma)}$ and $\h \in
G^{o(\tau)}$. Let $a \in o(\sigma)$, $b \in o(\tau)$ and $c \in
o(\sigma\tau)$. The multiplication of the group $G^I \rtimes
\Sigma_I$ yields a map
\begin{equation}\label{surjmap}
 \cO_{\g} \sigma \times \cO_{\h}\tau
\rightarrow \bigsqcup_{\w}\cO_{\w}\sigma \tau
\end{equation}
where the disjoint union over $\w$ runs over the elements of $
G^{o(\sigma \tau)}$ such that $\prod_c \w_c=\prod_a \g_a \prod_b
\h_b$. There are $(G^I\times G^I)$-actions on $\cO_{\g} \sigma
\times \cO_{\h}\tau$ and $\bigsqcup_{\w}\cO_{\w}\sigma \tau $, and
the map (\ref{surjmap}) is equivariant with respect to these
actions. Furthermore, the action of $G^I\times G^I$ on
$\bigsqcup_{\w}\cO_{\w}\sigma \tau $ is transitive and, in
particular, all of the fibers have the same cardinality.
\end{lem}

\Pf Let $r:=|o(\sigma\tau)|$ and
$o(\sigma\tau)=\{c_1,c_2,\cdots,c_r\}$. There is the action of
$G^I\times G^I$ on $\cO_{\g} \sigma \times \cO_{\h}\tau$ by
componentwise conjugation and, by the map (\ref{surjmap}), it
induces an action of $G^I\times G^I$ on
$\bigsqcup_{\w}\cO_{\w}\sigma \tau$, $i.e.$ for all $f_1,f_2 \in
G^I$,
\[
w \sigma\tau \stackrel{(f_1,f_2)}{\longrightarrow} f_1^{-1}
f_1^{\sigma^{-1}} (f_2^{-1})^{\sigma^{-1}}
f_2^{(\sigma\tau)^{-1}}w\sigma\tau.
\]
The map (\ref{surjmap}) is equivariant with respect to these
actions. The action of the subgroup generated by all diagonal
elements $(g,g) \in G^I\times G^I$ on $\bigsqcup_{\w}\cO_{\w}\sigma
\tau$ coincides with the action of $G^I$ by conjugation on
$\bigsqcup_{\w}\cO_{\w}\sigma \tau$. Hence, to prove the
transitivity of the $(G^I\times G^I)$-action on
$\bigsqcup_{\w}\cO_{\w}\sigma \tau $, it suffices to show that, for
a given $\w$ such that $\prod_c \w_c=\prod_a \g_a \cdot \prod_b
\h_b$, there exists an $f \in G$ such that $\epsilon_{\g}\sigma\cdot
f\epsilon_{\h}\tau f^{-1}$ belongs to $\cO_{\w}\sigma \tau$.
However, such an $f$ is a solution to the following set of $r$
equations for the $\{f_i\}_{i\in I}$ where $k=1,\cdots,r$:
\begin{equation}\label{2-1}
\w_{c_k} = \prod_{i\in c_k}
(\epsilon_{\g})_i(\epsilon_{\h})_{\tau(i)} \cdot
f_{\tau(i)}f_i^{-1}.
\end{equation}
Let us call the equation associated with $c_k$ the \emph{$k$-th
equation}. Let
$B_k:=(c_r\backslash\tau(c_k))\cup(\tau(c_k)\backslash c_k)$ and
then the $k$-th equation is an equation for $\{f_i\}_{i\in B_k}$.
Observe that the product of all $r$ equations is
$\prod_c\w_c=\prod_a\g_a\cdot\prod_b\h_b$.
Hence, if $f$ in $G^I$ satisfies the first $r-1$
equations, then it satisfies the $r$-th equation trivially.

Let $m=1,\cdots,r-1$. Since $o(\sigma\tau,\tau)=\{I\}$, there exists
$j_m \in B_m$ such that $j_m$ is not contained in $B_k$ for all
$k=1,\cdots,m-1$. If we are given $f_i \in G$ for all $i \in I
\backslash \{j_m\}_{m=k}^{r-1}$, the $k$-th equation determines
$f_{j_k}$ uniquely. Hence, once we choose $f_i$ in $G$ for all $i
\in I\backslash\{j_m\}_{m=1}^{r-1}$, by induction on $k$, we
uniquely find $\{f_{j_m}\}_{m=1}^{r-1}$ satisfying the set of
equations (\ref{2-1}).

Thus, the action of $G^I\times G^I$ on $\bigsqcup_{\w}\cO_{\w}\sigma
\tau $ is transitive and, in particular, the cardinality of each
fibre is $|G|^{n+1-|o(\sigma)|-|o(\tau)|}$. \qed

\section{{\bf The Lehn and Sorger algebra}}\label{LS}
In this section, we review the algebra $\cA\{\Sigma_I\}$ introduced by Lehn
and Sorger \cite{L-S} associated to a Frobenius algebra
$\cA$. In particular, $\cA$ could be the ordinary cohomology ring of
a compact almost complex manifold of complex dimension $d$. In this
paper, we will be primarily interested in the case where $\cA$ is
the orbifold cohomology ring of a global quotient of a compact
almost complex manifold of complex dimension $d$ by a finite group.

\begin{defn}\label{comultiplication}
Let $\cA$ be a Frobenius algebra. The associative multiplication
$\mu$ defines the \emph{multi-product} $ \cA^{\otimes n} \rightarrow
\cA$ by
\[
\bx_1\otimes\cdots\otimes\bx_n \mapsto \bx_1 \cdot \bx_2 \cdot\
\cdots \ \cdot \bx_n
\]
which will also be denoted by $\mu$. Let $\mu^{*}:\cA^{*}
\rightarrow \cA^{*}\otimes \cdots \otimes \cA^{*}$ be the dual of
the multi-product $\mu$. Let $\eta^{\flat}$ be a map which sends
elements of $\cA$ to $\cA^{*}$ by
\[
\bx \mapsto \eta(\ \  ,\bx).
\]
We can extend $\eta^{\flat}$ to the map $\cA
\otimes\cdots\otimes\cA\rightarrow\cA^{*}\otimes\cdots
\otimes\cA^{*}$ by applying $\eta^{\flat}$ to each factor, which we
also denote by $\eta^{\flat}$. We can define the
\emph{comultiplication} $\mathbf{m}_{*}:\cA \rightarrow \cA^{\otimes
n}$ by
\[
\cA \stackrel{\eta^{\flat}}{\longrightarrow} \cA^{*}
\stackrel{\mu^{*}}{\longrightarrow}
\cA^{*}\otimes\cdots\otimes\cA^{*}\stackrel{(\eta^{\flat})^{-1}}{\longrightarrow}
\cA\otimes\cdots \otimes\cA.
\]
\end{defn}

\begin{defn}\label{33}
Let $[n]:=\{1,\cdots,n\}$ and let $J$ be a finite set of cardinality
$n$. Let $\{A_i\}_{i \in J}$ be the collection of copies of $\cA$
indexed by $J$. $\cA^{\otimes J}$ is defined by
\[
\cA^{\otimes J}:= \left.\left( \bigoplus_{f:[n]\iso J}
\cA_{f(1)}\otimes \cdots \otimes \cA_{f(n)} \right)\right/ \Sigma_n
\]
where the direct sum runs over the set of all bijections $[n] \iso
I$ and the action of the permutation group $\Sigma_n$ of $[n]$ on
$\left( \bigoplus_{f:[n]\iso J} A_{f(1)}\otimes \cdots \otimes
A_{f(n)} \right)$ is induced by the bijections from $[n]$ to $J$.

For any finite sets $J_1$,$J_2$ and a surjective map
$\phi:J_1\rightarrow J_2$, the multi-product and comultiplication in
Definition \ref{comultiplication} can be generalized to the maps
$\phi^{*}:\cA^{J_1} \rightarrow \cA^{J_2}$ and $\phi_{*}:\cA^{J_2}
\rightarrow \cA^{J_1}$ respectively. Let $n:=|J_1|$  and $k:=|J_2|$.
Consider the ring homomorphism $\phi_{n,k} : \cA^{\otimes
n}\rightarrow \cA^{\otimes k}$ which sends $a_1\otimes \cdots
\otimes a_n$ to $(a_1\otimes\cdots\otimes a_{n_1})\otimes \cdots
\otimes (a_{n_1+\cdots +n_{k-1}+1} \otimes \cdots \otimes a_n)$
where $n_i:=|\phi^{-1}(g(i))|$.  Choose a bijection $g:[k]\iso J_2$
and then there is a bijection $f : [n]\iso J_1$ such that, for each
$i\in [k]$, one has $\phi^{-1}(g(i)) = \{ f(n_1 + \cdots + n_{i-1} +
1),\cdots , f(n_1 + \cdots + n_i)\}$. The composition
\[
\phi^{*}:\cA^{J_1} \iso \cA^{\otimes n}
\stackrel{\phi_{n,k}}{\longrightarrow} \cA^{\otimes k} \iso
\cA^{J_2}
\]
is independent of the choices of $f$ and $g$, where
the first and third maps are obvious isomorphisms induced by $f$ and
$g$. Let $\phi_{*}:\cA^{J_2} \rightarrow \cA^{J_1}$ be the linear map
adjoint to $\phi^*$ with respect to the metric induced from the
metric on $\cA$. In particular, if $J_1=[n]$ and $k=1$,
then $\phi^{*}$ and $\phi_{*}$ are the multi-product $\mu$ and
the comultiplication $\bm_{*}$ in Definition \ref{comultiplication}.
\end{defn}
\begin{defn}\label{Eulerclass}
The \emph{Euler class} $\e$ of $\cA$ is the image of $\1$ under the
composition of the maps
\[
\cA \stackrel{\bm_*}{\longrightarrow}
\cA\otimes\cA\stackrel{\mu}{\longrightarrow}\cA\ .
\]
\end{defn}

Now we define the algebra $\cA\{\Sigma_I\}$, following \cite{L-S}.
Let $\cA\{\Sigma_I\}$ be the $\Sigma_I$-graded vector space
\[
\cA\{\Sigma_I\} := \bigoplus_{\sigma \in \Sigma_I}\cA^{\otimes
o(\sigma)} \cdot \sigma \ .
\]
Let $\deg \bx$ be the $\Q$-grading of $\bx \in \cA^{\otimes
o(\sigma)}$ and let $d$ be the half of the top degree of $\cA$ as in
Definition \ref{Q-grading}. The $\Q$-grading of $\bx\sigma
\in\cA^{\otimes o(\sigma)} \cdot \sigma $ is defined by
\begin{equation}\label{LSdegree}
\deg^{LS}\bx\sigma:=(\deg \bx) + \frac{d \cdot l_{\sigma}}{2}
\end{equation}
where $l_{\sigma}$ is the length of the permutation $\sigma$.

We have a $\Sigma_I$-action $\rho$ on $\cA\{\Sigma_I\}$ which
preserves the $\Q$-grading, namely, the action of $\sigma \in
\Sigma_I$ on $I$ induces a bijection $o(\tau) \iso
o(\sigma^{-1}\tau\sigma)$ for each $\tau\in \Sigma_I$ and hence an
automorphism of $\cA\{\Sigma_I\}$,
 \[
\rho(\sigma):\cA^{o(\tau)}\cdot\tau\stackrel{\simeq}{\longrightarrow}
\cA^{o(\sigma^{-1}\tau\sigma)}\cdot \sigma^{-1}\tau\sigma.
\]
In particular, since the action of $\sigma$ on $o(\sigma)$ is the
identity map, the induced action $\rho(\sigma)$ on $\cA^{\otimes
o(\sigma)} \cdot \sigma$ is the identity map.

For all $\sigma,\tau\in \Sigma_I$, we have the following maps by
Definition \ref{33},
\[
f_1^{*}:\cA^{o(\sigma)} \rightarrow\cA^{o(\sigma,\tau)},\
f_2^{*}:\cA^{o(\tau)}\rightarrow \cA^{o(\sigma,\tau)} \ \text{ and
}\ \ f_{3*}: \cA^{o(\sigma,\tau)}\rightarrow\cA^{o(\sigma\tau)}
\]
induced respectively from the canonical surjections
\[
f_1:o(\sigma)\twoheadrightarrow o(\sigma,\tau),\
f_2:o(\tau)\twoheadrightarrow o(\sigma,\tau) \text{ and }
f_3:o(\sigma\tau) \twoheadrightarrow o(\sigma,\tau).
\]
The product on $\cA\{\Sigma_I\}$ is defined by
\begin{equation}\label{def LSprod}
\bx\sigma \cdot \by\tau := f_{3*}(f_1^{*}(\bx)\cdot
f_2^{*}(\by)\cdot \e^{gd(\sigma,\tau)})\cdot \sigma\tau
\end{equation}
where
\[
\e^{gd(\sigma,\tau)}:=\bigotimes_{c \in o(\sigma,\tau)}
\e^{gd(\sigma,\tau)_c}
\]
and
\begin{equation}\label{graph defect}
gd(\sigma,\tau)_c:=\frac{1}{2}(|c|+2-|c/\lan\sigma\ran
|-|c/\lan\tau\ran|-|c/\lan\sigma\tau\ran|).
\end{equation}
$gd(\sigma,\tau)_c$ is called the \emph{graph defect} of $\sigma$
and $\tau$ on $c \in o(\sigma, \tau)$ and is a non-negative integer
by Lemma 2.7 in \cite{L-S}. By Proposition 2.13 in \cite{L-S}, the
multiplication defined by Equation (\ref{def LSprod}) is associative
and $\Sigma_I$-equivariant. By Proposition 2.14 in \cite{L-S}, the
multiplication is also braided commutative. The metric $\eta$ is
defined by
\begin{equation}
\eta(\bx\sigma, \by\tau):=
\begin{cases}
\eta(\1, \bx\sigma\cdot\by\tau) & \text{ if } \sigma\tau=id \\
0 & \text{ otherwise,}
\end{cases}
\end{equation}
where $\eta$ on the right-hand side of the equality is induced from
the metric on $\cA^I$. This metric is non-degenerate and
$\Sigma_I$-invariant by Proposition 2.16 in \cite{L-S}. The braided
commutativity and the self-invariance axioms imply that this metric is
symmetric. The associativity of the product implies the invariance
of the metric. Thus, we have shown the following.
\begin{prop}
$\cA\{\Sigma_I\}$ satisfies all the axioms of a $\Q$-graded
$\Sigma_I$-Frobenius algebra except, possibly, the trace axiom.
\end{prop}

\begin{rem}
If $\cA$ is connected, that is, the subspace of all elements
with trivial $\Q$-grading is $1$-dimensional, then
it is straightforward to prove that $\cA\{\Sigma_I\}$
satisfies the trace axiom.
Indeed, the traces of $L_v\circ\rho(\tau^{-1})$ and
$\rho(\sigma)\circ L_v$ are zero
for any homogenous element $v\in\cA^{o([\sigma,\tau])}\cdot[\sigma,\tau]$
unless $v=1$ and $[\sigma,\tau]=id$.
Hence, the trace axiom is trivially satisfied.
\end{rem}

\begin{rem}
We will later prove that $\cA\{\Sigma_I\}$ satisfies
the trace axiom if $\cA$ is the orbifold cohomology of a global
quotient of a compact almost complex manifold with an action of a
finite Abelian group, or if $\cA$ is the center of the group ring of
any finite group.
\end{rem}

\begin{rem}\label{splitting LS}
Let $\lambda:=\{d\}$ be a partition of $I$, $i.e.$ for all
$d\not=d'\in \lambda$, $d\cap d'$ is empty and $\cup_{d \in \lambda}
d=I$. A partition $\lambda'$ of $I$ is a \emph{subpartition} of
$\lambda$, denoted by $\lambda' < \lambda$, if and only if each $d'
\in \lambda'$ is contained in some $d \in \lambda$. Consider the
following subspace of $\cA\{\Sigma_I\}$
\[
\cA\{\Sigma_I\}(\lambda):=\bigoplus_{o(\sigma)<\lambda}
\cA^{o(\sigma)}\cdot\sigma.
\]
It is clear that this subspace is actually a subalgebra and we can
show that there is a ring isomorphism,
\begin{equation}\label{splitting LS product}
\cA\{\Sigma_I\}(\lambda) \cong
\bigotimes_{d\in\lambda}\cA\{\Sigma_d\}.
\end{equation}
In fact, if $o(\sigma)<\lambda$ and $o(\tau)<\lambda$,
then we have $o(\sigma,\tau)<\lambda$. Let $f_{\lambda}
:o(\sigma,\tau)\twoheadrightarrow \lambda$ be the obvious
surjection. Let $r:=|\lambda|$ and choose a bijection $f:[r] \iso
\lambda$. Let $c_k:=f(k)$. If $\bx=\bx_1\otimes\cdots\otimes \bx_r$ and
$y=\by_1\otimes\cdots\otimes\by_r$ for
$\bx_k\in\cA^{(f_{\lambda}\circ f_1)^{-1}(c_k)}$ and
$\by_k\in\cA^{(f_{\lambda}\circ f_2)^{-1}(c_k)}$, then we have
\begin{equation}
\bx\sigma \cdot \by\tau = \bigotimes_{k=1}^{r}
f_{3,k*}(f_{1,k}^{*}(\bx_k)\cdot f_{2,k}^{*}(\by_k)\cdot
\e^{gd(\sigma,\tau)_{c_k}})\cdot \sigma\tau
\end{equation}
where $f_{i,k}:(f_{\lambda}\circ f_i)^{-1}(c_k)\twoheadrightarrow
\{k\}$ for $i=1,2,3$ are the obvious surjections induced by $f$,
$f_{\lambda}$ and $f_i$'s.
\end{rem}

\section{{\bf The Lehn and Sorger
algebra associated to orbifold cohomology}}\label{Orbifold}

We review the definition of the stringy cohomology and the orbifold
cohomology, following \cite{FG} and \cite{JKK2}. Let $X$ be a
compact almost complex manifold with an action of a finite group
$G$ preserving the almost complex structure. Denote the action of
$G$ on $X$ by $\rho$. For any set of $r$ elements in $G$,
$\{\alpha_1,\alpha_2, \cdots, \alpha_r\}$, we denote the fixed point locus of $\lan
\alpha_1,\cdots,\alpha_r\ran $ in $X$ by $X^{\alpha_1,\cdots,\alpha_r}$. Define the
\emph{inertia manifold} of $X$ by
\[
I_G(X):=\{(x,\alpha) \in X\times G \ |\ \rho(\alpha)x=x\}
=\bigsqcup_{\alpha \in G} X^{\alpha}.
\]
Let $\cH(X,G)$ be the ordinary cohomology $H^{*}(I_G(X))$ of
$I_G(X)$ and it is a $G$-graded $G$-module where
\[
\cH(X,G)=\bigoplus_{\alpha\in G}\cH^X_{\alpha}
:=\bigoplus_{\alpha\in G} H^{*}(X^{\alpha}).
\]
The subspace generated by vectors that are graded by non-trivial
group elements, $\bigoplus_{\alpha\not=1}\cH^X_{\alpha}$, is called
the \emph{twisted sector}. Let $\alpha,\beta \in G$ and let
$\bq:X^{\alpha,\beta}\inc X^{\alpha\beta}$ be the canonical
inclusion map. For $\bx \in H^{*}(X^{\alpha})$ and $\by \in
H^{*}(X^{\beta})$, their product is defined by
\begin{equation}\label{product}
\bx\cdot \by := \bq_{*}\left[\bx|_{X^{\alpha,\beta}}\cup
\by|_{X^{\alpha,\beta}}\cup c_{top}\left(\sR(\alpha,\beta)\right)
\right]
\end{equation}
where $\sR(\alpha,\beta)$ is the obstruction bundle over
$X^{\alpha,\beta}$ introduced in \cite{FG}. Let
\begin{equation}\label{virtual class}
c(\alpha,\beta):=c_{top}\left(\sR(\alpha,\beta)\right).
\end{equation}
In \cite{JKK2}, the following equality in the $K$-theory of
$X^{\alpha,\beta}$, $K(X^{\alpha,\beta})$, is shown:
\begin{equation}\label{obstruction bundle}
\sR(\alpha,\beta)=TX^{\alpha,\beta}\ominus TX|_{X^{\alpha,\beta}}
\oplus \sS_{\alpha}|_{X^{\alpha,\beta}} \oplus
\sS_{\beta}|_{X^{\alpha,\beta}}
\oplus\sS_{(\alpha\beta)^{-1}}|_{X^{\alpha,\beta}}.
\end{equation}
Here, the class $\sS_{\alpha}$ in $K(X^{\alpha})$ is defined by

\begin{equation}\label{def-Sbundle}
\sS_{\alpha}:=\bigoplus_{k=0}^{r-1} \frac{k}{r}W_{\alpha,k}
\end{equation}
where $r$ is the order of $\alpha$, and $W_{\alpha,k}$ is the
eigenbundle of $W_{\alpha}:=TX|_{X^{\alpha}}$ such that $\alpha$
acts with the eigenvalue $\exp(-2\pi ki/r)$. In particular, the rank
of $\sS_{\alpha}$ on a connected component $C$ of $X^{\alpha}$ is
called \emph{age} of $\alpha$ on $C$ and is denoted by
$age(\alpha)_C$. It is worth noting that, for every $\alpha,\beta$
and $\gamma \in G$,
\begin{equation}\label{equiv S}
\sS_{\alpha}=\rho(\beta)^{*}\sS_{\beta^{-1}\alpha\beta}
\end{equation}
where $\rho(\beta):X^{\alpha} \stackrel{\simeq}{\longrightarrow}
X^{\beta^{-1}\alpha\beta}$, and
\begin{equation}\label{G-equivR}
\rho(\gamma)_{\ast}\sR(\alpha,\beta)=\sR\left(\gamma^{-1}
\alpha\gamma,\gamma^{-1}\beta\gamma\right)
\end{equation}
where $\rho(\gamma):X^{\alpha,\beta}\iso
X^{\gamma^{-1}\alpha\gamma,\gamma^{-1}\beta\gamma}$.

 The metric $\eta$ of $\cH(X,G)$ is defined by
\[
\eta(\bx,\by):=\int_{I_G(X)} \bx \cup \iota^{*}\by
\]
where $\iota: I_G(X) \rightarrow I_G(X)$ is the canonical
$G$-equivariant involution on $I_G(X)$ taking $X^{\alpha}$ to
$X^{\alpha^{-1}}$. With the product and the metric above, $\cH(X,G)$
becomes a $G$-Frobenius algebra and is called the \emph{stringy
cohomology} of $G$-space $X$. For a global quotient $[X/G]$, the
$G$-coinvariants of the stringy cohomology is isomorphic as a
Frobenius algebra to the orbifold cohomology of Chen-Ruan \cite{CR},
$i.e.$
\[
 \cH(X,G)^{G} = H_{orb}^{*}([X/G]).
\]

\begin{lem}
Let $X$ be a compact almost complex manifold with an action of a
finite group $G$. Let $\e$ be the Euler class, as defined in
Definition \ref{Eulerclass}, of the orbifold cohomology
$H_{orb}^{*}([X/G])$. We have
\[
\e=\sum_{\alpha,\beta \in G, \alpha\beta=\beta\alpha}
\br_{\alpha,\beta*}\ c_{top}(TX^{\alpha,\beta})
\]
where $\br_{\alpha,\beta} : X^{\alpha,\beta}\inc X$ is the canonical
inclusion for all $\alpha,\beta$ in $G$.
\end{lem}

\Pf Let $X_1/G,\cdots, X_r/G$ be connected components of the
quotient $X/G$ and then each of them is itself an orbifold. By Definition
\ref{Eulerclass}, the Euler class of $H^*_{orb}([X/G])$ is
$\e=\e_1+\cdots+\e_r$ where $\e_k$ is the Euler class of
$H_{orb}^*([X_k/G])$. Hence we can assume that the quotient space
$X/G$ is connected without loss of generality.

Let $\chi(M/\Gamma)$ be the $\Gamma$-equivariant Euler characteristic for a
compact manifold $M$ with an action of a finite group $\Gamma$. The
following identity is well-known (c.f. \cite{AS}):
\begin{equation}\label{finite-Euler}
\chi(M/\Gamma)=\frac{1}{|\Gamma|} \sum_{\alpha\in \Gamma}
\chi(M^\alpha).
\end{equation}

Let $C_{\alpha}$ be the conjugacy class of $\alpha$ in $G$ and
$\mathbf{vol}$ is the $G$-invariant class of a volume form of $X$.
By Definition \ref{Eulerclass},
\[
\e=|G| \sum_{\alpha \in G} \frac{1}{|C_{\alpha}|}
\chi(X^{\alpha}/\cZ_G(\alpha))\mathbf{vol}.
\]
By Equation (\ref{finite-Euler}), we obtain
\begin{equation*}
\e=\sum_{\alpha\beta=\beta\alpha} \chi(X^{\alpha,\beta})\mathbf{vol}
= \sum_{\alpha\beta=\beta\alpha} \br_{\alpha,\beta*}
c_{top}(TX^{\alpha,\beta}).
\end{equation*}
\qed

We compute the multi-product $\mu:H^{*}_{orb}([X/G])^{\otimes n}
\rightarrow H^{*}_{orb}([X/G])$ and comultiplication
$\mathbf{m}_{*}:H^{*}_{orb}([X/G])\rightarrow
H^{*}_{orb}([X/G])^{\otimes n}$ given in Definition
\ref{comultiplication} in the next two propositions.

\begin{rem}\label{notation kunneth}
For all $g \in G^n$, let $\Delta_n : X^{g_1,\cdots, g_n} \rightarrow
X^{g_1} \times \cdots \times X^{g_n}$ be the diagonal embedding. Let
$\bx_k \in H^{*}(X^{g_k})$ for all $k=1,\cdots, n$. We regard $\bx_1
\otimes \cdots \otimes \bx_n$ as belonging to $H^{*}((X^n)^g)$ by
the K\"{u}nneth theorem. We have
\[
\bx_1|_{Z_n} \cup \cdots \cup \bx_n|_{Z_n}=\Delta_n^{*}(\bx_1
\otimes \cdots \otimes \bx_n)
\]
where $Z_n:=X^{g_1,\cdots, g_n}$.
\end{rem}

\begin{prop}(Multi-product)\label{multiprod}
Suppose that $G$ is an Abelian group. For all $g \in G^n$ and
$\bx_1\otimes\cdots\otimes\bx_n\in H^{*}((X^n)^{g})^{G^n}$, we have
\begin{equation}\label{multiproductformula}
\mu(\bx_1 \otimes \cdots \otimes \bx_n) =\br_{n*}
\left(\Delta_n^{*}(\bx_1 \otimes \cdots \otimes \bx_n)\cup c(g_1,
\cdots, g_n) \right)
\end{equation}
where $\br_n:Z_n \inc X^{g_1\cdots g_n}$ is the canonical inclusion
and $c(g_1, \cdots, g_n)$ is the top Chern class of the vector
bundle which is equal to the following element in $K(Z_n)$:
\begin{equation}\label{MultiObstructionBundle}
TZ_n \ominus TX|_{Z_n} \oplus  \bigoplus_{k=1}^n\sS_{g_i}|_{Z_n}
\oplus \sS_{(g_1\cdots g_n)^{-1}}|_{Z_n}.
\end{equation}
\end{prop}

\Pf We will prove the proposition by induction on $n$. When $n=1$,
the Equation (\ref{multiproductformula}) is trivial. Let $g_k \in G$
and $\bx_k \in H^{*}(X^{g_k})$ for all $k=1,\cdots, n$ and suppose
that
\begin{equation}\label{47}
\mu(\bx_1 \otimes \cdots \otimes \bx_{n-1}) =\br_{n-1*} \left(
\Delta_{n-1}^{*}(\bx_1 \otimes \cdots \otimes \bx_{n-1}) \cup
c(g_1,\cdots, g_{n-1}) \right).
\end{equation}
Consider the following commuting diagram
\begin{equation}\label{excess}
\begin{CD}
Z_{n-1} @>\br_{n-1}>> X^{g_1\cdots g_{n-1}} \\
@A\iota_1 AA   @AA\iota_3A  \\
Z_n  @ >>\iota_2> X^{g_1\cdots g_{n-1}, g_n}
\end{CD}
\end{equation}
and the maps
\begin{equation}
X^{g_n} \stackrel{\gamma_1}{\longleftarrow} X^{g_1\cdots
g_{n-1}, g_n} \stackrel{\gamma_2}{\longrightarrow} X^{g_1\cdots g_n}
\end{equation}
where all of the maps are the obvious inclusions. The excess
intersection formula \cite{Qu} associated to diagram
(\ref{excess}) yields the identity in $H^{*}(Z_n)$,
\begin{equation}\label{excess formula}
\iota_3^{*}\br_{n-1*}(\bx)=\iota_{2*}( \iota_1^{*} ( \bx) \cup E_n)
\end{equation}
where $E_n= c_{top}\left( TX^{g_1\cdots g_{n-1}}|_{Z_n} \oplus TZ_n
\ominus TX^{g_1\cdots g_{n-1},g_n}|_{Z_n} \ominus TZ_{n-1}|_{Z_n}
\right)$. By associativity of the product, we can write
$\mu(\bx_1\otimes\cdots\otimes\bx_n)=\mu(\bx_1\otimes \cdots\otimes
\bx_{n-1}) \cdot \mu(\bx_n)$. The right-hand side is computed as
follows.
\begin{eqnarray*}
&&\mu(\bx_1\otimes \cdots\otimes \bx_{n-1}) \cdot \mu(\bx_n) \\
&=&\gamma_{2*}\left[\iota_3^{*}
\br_{n-1*}\left(\Delta_{n-1}^{*}(\bx_1 \otimes \cdots \otimes
\bx_{n-1}) \cup c(g_1,\cdots, g_{n-1}) \right) \cup \gamma_1^{*}
\bx_{m} \cup c(g_1\cdots g_{n-1}, g_n)\right] \\
&=& \gamma_{2*}\left[\iota_{2*}\left(  \iota_1^{*} \left(
\Delta_{n-1}^{*}(\bx_1 \otimes \cdots \otimes \bx_{n-1})\cup
c(g_1,\cdots, g_{n-1}) \right) \cup E_n \right) \cup \gamma_1^{*}
\bx_{n} \cup c(g_1\cdots g_{n-1}, g_n)
 \right] \\
&=&(\gamma_2\iota_2)_{*}\left[\iota_1^{*}\left(
\Delta_{n-1}^{*}(\bx_1 \otimes \cdots \otimes \bx_{n-1})\cup
c(g_1,\cdots, g_{n-1}) \right)\cup E_n\cup\iota_2^{*}(\gamma_1^{*}
\bx_{n}
\cup c(g_1\cdots g_{n-1},g_n))\right] \\
&=&\br_{n*}\left[\Delta_n^{*}(\bx_1 \otimes \cdots \otimes
\bx_n)\cup c(g_1,\cdots, g_{n-1}) |_{Z_n}\cup E_n\cup c(g_1\cdots
g_{n-1}, g_n) |_{Z_n}\right]
\end{eqnarray*}
where the first equality follows from the definition of the product
in Equation (\ref{product}) and the induction hypothesis, the second
equality follows from Equation (\ref{excess formula}) and the third
follows from the projection formula. Finally,
\[
c(g_1,\cdots, g_{n-1})|_{Z_n}\cup E_n\cup c(g_1\cdots
g_{n-1},g_n)|_{Z_n}
\]
is equal to the top Chern class of the bundle which belongs to the
following class in $K(Z_n)$,
\begin{eqnarray*}
&& TZ_{n-1}|_{Z_n} \ominus TX|_{Z_n} \oplus
\bigoplus_{k=1}^{n-1}\sS_{g_i}|_{Z_n} \oplus \sS_{(g_1\cdots
g_{n-1})^{-1}}|_{Z_n}\\
&& \oplus TX^{g_1\cdots g_{n-1}}|_{Z_n} \oplus TZ_n \ominus
TX^{g_1\cdots g_{n-1},g_n}|_{Z_n} \ominus
TZ_{n-1}|_{Z_n}\\
&& \oplus   TX^{g_1\cdots g_{n-1},g_n}|_{Z_n} \ominus  TX|_{Z_n}
\oplus \sS_{g_1\cdots g_{n-1}}|_{Z_n} \oplus \sS_{g_n}|_{Z_n} \oplus
\sS_{(g_1\cdots g_n)^{-1}}|_{Z_n}.
\end{eqnarray*}
This class in $K(Z_n)$ simplifies to
\[
TZ_n \ominus TX|_{Z_n} \oplus  \bigoplus_{k=1}^n\sS_{g_i}|_{Z_n}
\oplus \sS_{(g_1\cdots g_n)^{-1}}|_{Z_n}
\]
where we used the identity
\[
\sS_{g_1\cdots g_{n-1}}\oplus \iota^{*}\sS_{(g_1\cdots g_{n-1})^{-1}
}= TX|_{X^{g_1\cdots g_{n-1}}} \ominus TX^{g_1\cdots g_{n-1}}.
\]
 \qed

\begin{prop}(Co-multiplication)\label{coprod}
Suppose that $G$ is an Abelian group. For all $\bx_{\alpha} \in
H^{*}(X^{\alpha})^G$, we have
\[
\bm_{*}(\bx_{\alpha})=\frac{1}{|G|} \sum_{f \in G^n}
\sum_{g \in G^n} \rho(g)[\Delta_{*}(\br^{*} \bx_{\alpha} \cup
c(f_1^{-1}, \cdots,f_n^{-1}))]
\]
where the first sum runs over all elements $f$ in $G^n$ such that
$f_1\cdots f_n=\alpha$, and where $\br:X^{f_1,\cdots, f_n} \inc
X^{\alpha}$ is the canonical inclusion and $\Delta:X^{f_1,\cdots,
f_n} \rightarrow X^{f_1^{-1}} \times \cdots \times X^{f_n^{-1}}$ is
the diagonal embedding.
\end{prop}

\Pf The comultiplication $\bm_{*}$ restricted to
$H^{*}(X^{\alpha})^G$ is defined by the following commuting diagram:
\begin{equation}
\begin{CD}
H^{*}(X^\alpha)^G @>\bm_{*}>> \bigoplus_{f \in G^n}
H^{*}(X^{f_1} \times \cdots \times X^{f_n})^{G^n} \\
@V\eta^{\flat}VV @AA(\eta^{\flat})^{-1}A \\
\left(H^{*}(X^{\alpha^{-1}})^G\right)^{*} @>\mu^{*}>>\bigoplus_{f
\in G^n} \left(H^{*}(X^{f_1^{-1}} \times \cdots \times
X^{f_n^{-1}})^{G^n}\right)^{*}.
\end{CD}
\end{equation}
Let $\pi : X^{\alpha^{-1}} \rightarrow \{pt\}$ and $\pi':
X^{f_1}\times \cdots \times X^{f_n} \rightarrow \{pt\}$ be the
obvious projection maps. For all $\bx_{\alpha} \in
H^{*}(X^{\alpha})^G$, $\eta^{\flat}(\bx_{\alpha})$ is the linear
functional on $H^{*}(X^{\alpha^{-1}})^G$ taking
\[
 \bx_{\alpha^{-1}} \mapsto
\frac{1}{|G|}\pi_{*}(\bx_{\alpha^{-1}} \cup \iota^{*}\bx_{\alpha}).
\]
Since $\mu^{*}$ is the dual of the multiplication, $\mu^{*}\circ
\eta^{\flat}(\bx_{\alpha})$ is the linear functional on
$H^{*}_{orb}([X/G])^{\otimes n}$ which sends $\by \in
H^{*}_{orb}([X/G])^{\otimes n}$ to $\frac{1}{|G|}\pi_{*}(\mu(\by)
\cup \iota^{*}\bx_{\alpha})$. Applying Proposition \ref{multiprod},
we can write
\[
\frac{1}{|G|}\pi_{*}(\mu(\by) \cup \iota^{*}\bx_{\alpha}) =
\sum_{f\in G^n}\frac{1}{|G|}\pi_{*} \left(\br_{*} \left(\Delta^{*}
\by \cup c(f_1^{-1},\cdots,f_n^{-1})\right)\cup
\iota^{*}\bx_{\alpha}\right)
\]
where the sum on the right-hand side runs over $f \in G^n$ such that
$f_1\cdots f_n=\alpha$. Furthermore, using the projection formula,
the right-hand side is equal to
\begin{eqnarray*}
&& \sum_{f \in G^n}\frac{1}{|G|}\pi_{*}\circ \br_{*}
\left(\Delta^{*}\by \cup c(f_1^{-1},\cdots,f_n^{-1})\cup \br^{*}\bx_{\alpha}\right) \\
&=& \sum_{f \in G^n}\frac{1}{|G|}\pi'_{*}\circ \Delta_{*}
\left(\Delta^{*}\by \cup c(f_1^{-1},\cdots,f_n^{-1}) \cup \br^{*}\bx_{\alpha}\right) \\
&=& \sum_{f \in G^n}\frac{1}{|G|}\pi'_{*} \left(\by \cup\Delta_{*}
\left(c(f_1^{-1},\cdots,f_n^{-1}) \cup \br^{*}\bx_{\alpha}\right)\right) \\
&=& \sum_{f \in G^n}\frac{1}{|G|}\pi'_{*} \left(\frac{1}{|G|^n}\sum_{g \in G^n}\rho(g)(\by)
\cup\Delta_{*}\left(c(f_1^{-1},\cdots,f_n^{-1}) \cup \br^{*}\bx_{\alpha}\right)\right)\\
&=& \frac{1}{|G|^n}\pi'_{*} \left(\by
\cup\frac{1}{|G|}\sum_{f\in G^n}\sum_{g \in G^n} \rho(g) \left[\Delta_{*}
\left(c(f_1^{-1},\cdots,f_n^{-1}) \cup \br^{*}\bx_{\alpha}\right)\right]\right) \\
\end{eqnarray*}
where the first equality follows from the commutative diagram
\begin{equation*}
\begin{CD}
X^{f_1,\cdots, f_n} @>\Delta>> X^{f_1} \times \cdots \times X^{f_n} \\
@V\br VV @ V\pi' VV \\
X^{\alpha} @>\pi >> \{pt\}.
\end{CD}
\end{equation*}
The second equality is obtained by the projection formula and the
third equality follows from the invariance of the metric and the
fact that $\by$ is $G$-invariant. The fourth equality follows from the invariance
and the linearity of metric.
\qed
\begin{rem}\label{general coprod and comult}
Propositions \ref{multiprod} and \ref{coprod} can be generalized to
the maps $\phi^*$ and $\phi_*$ in Definition \ref{comultiplication}
where $\phi : J_1 \rightarrow J_2$ is a surjective map of sets. For
all $\g \in G^{J_1}$, define $\overline{\phi}(\g) \in G^{J_2}$ by
\[
\overline{\phi}(\g)_j:=\prod_{i \in \phi^{-1}(j)} \g_i
\]
for all $j \in J_2$. Let
\[
\br_{\g}:\prod_{j \in J_2} X^{\lan \g_i \, |\, i \in \phi^{-1}(j)
\ran} \longrightarrow \prod_{j \in J_2} X^{\overline{\phi}(\g)_j}
\]
be the canonical inclusion and let
\[
\Delta_{\g}:\prod_{j \in J_2} X^{\lan \g_i\, |\,i \in \phi^{-1}(j)
\ran} \rightarrow (X^{J_1})^{\g}
\]
where $\Delta_{\g}$ restricted on $X^{\lan \g_i \,|\,i \in
\phi^{-1}(j) \ran}$ is the diagonal map to $\prod_{i \in
\phi^{-1}(\g)} X^{\g_i}$.

The multi-product $\mu : H^{*}((X^{J_1})^{\g})^{G^{J_1}} \rightarrow
H^{*}((X^{J_2})^{\phi(\g)})^{G^{J_2}}$ can be written as
\[
\mu(\bx)=\br_{\g*}\left(\Delta_{\g}^{*}(\bx)\cup \bigotimes_{j \in
J_2} c\left(\{\g_i\}_{i \in \phi^{-1}(j)}\right)\right).
\]
The comultiplication $\bm_{\ast}:H^{*}((X^{J_2})^{\h})^{G^{J_2}}
\rightarrow H^{*}_{orb}([X/G])^{\otimes J_1}$ can be written as
\[
\bm_{*}(\bx)=\frac{1}{|G|} \sum_{\f \in G^{J_1}}
 \sum_{\g \in G^{J_1}}
\rho(\g)\left[\Delta_{\f^{-1}*}\left(\br_{\f}^{*} \bx \cup
\bigotimes_{j \in J_1}c\left(\{\f_i^{-1}\}_{i \in
\phi^{-1}(j)}\right)\right)\right]
\]
where the first sum runs over all of the elements $\f \in G^{J_1}$
such that $\phi(\f)=\h$.
\end{rem}

We henceforward adopt the following notation:
\begin{notation}\label{notation fix locus}
Let $I_1, I_2, \cdots, I_n$ be finite sets. Let $\g_{I_k} \in
G^{I_k}$ where $k=1,2, \cdots, n$. Let $\lan \g_{I_1},\cdots
,\g_{I_n}\ran$ be the subgroup of $G$ generated by all of the
components of the $\g_{I_k}$'s. We denote the fixed point locus of
$\lan \g_{I_1},\cdots ,\g_{I_n}\ran$ by
$X^{\g_{I_1},\g_{I_2},\cdots,\g_{I_n}}$. Let $J$ be a finite set.
For all $\g \in G^{J}$, define $\sS_{\g} := \bigoplus_{j\in J}
\sS_{\g_j}|_{X^{\g}}$. If $G$ is Abelian, define
$c(\g_{I_1},\g_{I_2},\cdots,\g_{I_n})$ to be the top Chern class of
the vector bundle representing the element in $K(Z)$,
\[
TZ \ominus TX|_Z \oplus \bigoplus_{k=1}^n \sS_{\g_{I_k}}|_Z \oplus
\sS_{\alpha^{-1}}|_Z.
\]
Here $Z:=X^{\g_{I_1},\g_{I_2},\cdots,\g_{I_n}}$ and
$\alpha:=\prod_{k=1}^n\prod_{i \in I_k}\left(\g_{I_k}\right)_i \in
G$.

For example, if $\g:=(\g_1,\g_2)$ and $\h:=(\h_1,\h_2,\h_3)$, then
$X^{\g,\h}=X^{\g_1,\g_2,\h_1,\h_2,\h_3}$,
\[
\sS_{\g}=\sS_{\g_1}|_{X^{g_1,g_2}}\oplus \sS_{\g_2}|_{X^{g_1,g_2}}
\]
and
\[
c(\g,\h)=c_{top}\left(TZ \ominus TX|_Z \oplus \sS_{\g_1}|_Z \oplus
\sS_{\g_2}|_Z \oplus \sS_{\h_1}|_Z \oplus \sS_{\h_2}|_Z \oplus
\sS_{\h_3}|_Z \oplus \sS_{(\g_1\g_2\h_1\h_2\h_3)^{-1}}|_Z \right)
\]
where $Z:=X^{\g_1,\g_2,\h_1,\h_2,\h_3}$.
\end{notation}

\begin{lem}\label{LS prod geom discription}
Let $G$ be an Abelian group and let $\sigma, \tau \in \Sigma_I$.
Suppose that $\lan\sigma,\tau\ran$ acts on $I$ transitively and
let $gd:=gd(\sigma,\tau)_I$. Let $\g
\in G^{o(\sigma)}$ and $\h \in G^{o(\tau)}$. Let
$Z_{\w}:=X^{\g,\h,\w}$ for each $\w \in G^{o(\sigma\tau)}$ such that
$\prod_c \w_c =\prod_a\g_a \prod_b \h_b$. Let $\bq_{\w}$ and
$\bp_{\w}$ be the diagonal embeddings
\[
(X^{o(\sigma)})^{\g}\times(X^{o(\sigma)})^{\h}
\stackrel{\bp_{\w}}{\longleftarrow} Z_{\w}
\stackrel{\bq_{\w}}{\longrightarrow} (X^{o(\sigma\tau)})^{\w}.
\]
If $\bx\in H^{*}((X^{o(\sigma)})^{\g})^{G^{o(\sigma)}}$ and $\by\in
H^{*}((X^{o(\sigma)})^{\h})^{G^{o(\tau)}}$, then
\begin{equation}\label{61}
\bx \sigma\cdot \by \tau =\frac{1}{|G|} \sum_{\fa\in
G^{o(\sigma\tau)}} \rho(\fa) \left[ \sum_{\w} \bq_{\w*} \left(
\bp_{\w}^{*}(\bx \otimes \by) \cup c(\g,\h)|_{Z_{\w}} \cup
\e^{{gd}}|_{Z_{\w}} \cup c(\w^{-1})|_{Z_{\w}} \cup E_{\w} \right)
\right]\cdot \sigma\tau
\end{equation}
where the sum over $\w$ runs over all elements of
$G^{o(\sigma\tau)}$ such that $\prod_c \w_c =\prod_a\g_a \prod_b
\h_b$, and
\[
E_{\w}:=c_{top}\left (TX^{\prod_c \w_c}|_{Z_{\w}} \oplus TZ_{\w}
\ominus TX^{\g,\h}|_{Z_{\w}} \ominus X^{\w}|_{Z_{\w}}\right).
\]
\end{lem}

\Pf Consider the following diagram of the obvious inclusions:
\begin{equation}\label{excess2}
\begin{CD}
X^{\g,\h} @>\br_{\g,\h}>> X^{\prod_c \w_c} \\
@A\bs_1 AA   @AA\br_{\w}A  \\
Z_{\w}  @ >>\bs_2> X^{\w}.
\end{CD}
\end{equation}
The excess intersection formula associated to the above diagram
 yields the following identity in $H^{*}(X^{\w})$: for all $\alpha \in
 H^{*}(X^{\g,\h})$,
\begin{equation}\label{excess formula2}
\br_{\w}^{*}\circ \br_{\g,\h*}(\alpha) =
\bs_{2*}\left(\bs_1^{*}(\alpha)\cup E_{\w}\right).
\end{equation}
 Consider the following
sequence of maps,
\begin{equation}
\prod_a X^{\g_a} \times \prod_b X^{\h_b}
\stackrel{\Delta_{\g,\h}}{\longleftarrow} X^{\g,\h}
\stackrel{\br_{\g,\h}}{\longrightarrow}  X^{\prod_c \w_c}
\stackrel{\br_{\w}}{\longleftarrow}  X^{\w}
\stackrel{\Delta_{\w}}{\longrightarrow} \prod_c X^{\w_c}
\end{equation}
where $\Delta_{\g,\h}$ and $\Delta_{\w}$ are the diagonal
embeddings. By the definition of the product in the Lehn-Sorger algebra
and by Remark \ref{general coprod and comult} (see also Proposition
\ref{multiprod} and \ref{coprod}),
\begin{eqnarray*}
&&\bx \sigma\cdot \by \tau \\
&=& \frac{1}{|G|} \sum_{\fa\in G^{o(\sigma\tau)}} \rho(\fa) \left[
\sum_{\w} \Delta_{\w*} \left( \br_{\w}^{*} \left[ \br_{\g,\h *}
\left( \Delta_{\g,\h}^{*}(\bx\otimes\by) \cup c(\g,\h) \right)
\right] \cup \e^{{gd}}|_{X^{\w}} \cup c(\w^{-1}) \right)
\right]\\
&=& \frac{1}{|G|} \sum_{\fa\in G^{o(\sigma\tau)}} \rho(\fa)
\left[\sum_{\w} \Delta_{\w*} \left( \bs_{2*} \left[ \bs_1^{*} \left(
\Delta_{\g,\h}^{*}(\bx\otimes\by) \cup c(\g,\h) \right) \cup E_{\w}
\right] \cup \e^{{gd}}|_{X^{\w}} \cup c(\w^{-1}) \right)
\right]\\
&=& \frac{1}{|G|} \sum_{\fa\in G^{o(\sigma\tau)}} \rho(\fa)
\left[\sum_{\w} \Delta_{\w*} \circ \bs_{2*} \left( \bs_1^{*} \left(
\Delta_{\g,\h}^{*}(\bx\otimes\by) \cup c(\g,\h) \right) \cup E_{\w}
\cup \bs_1^{*}\e^{{gd}}|_{X^{\w}} \cup \bs_1^{*}c(\w^{-1}) \right)
\right]\\
&=& \frac{1}{|G|} \sum_{\fa\in G^{o(\sigma\tau)}} \rho(\fa)
\left[\sum_{\w} \bq_{\w*} \left( \bp_{\w}^{*}(\bx\otimes\by) \cup
c(\g,\h)|_{Z_{\w}} \cup E_{\w} \cup \e^{{gd}}|_{Z_{\w}}
\cup c(\w^{-1})|_{Z_{\w}} \right)
\right]\\
\end{eqnarray*}
where the second equality is obtained by the excess intersection
formula (\ref{excess formula2}), the third follows from the
projection formula and the fourth is obtained from the equalities
$\Delta_{\w}\circ \bs_2=\bq_{\w}$ and
$\bs_1\circ\Delta_{\g,\h}=\bp_{\w}$. \qed

\begin{rem}
By Remark \ref{splitting LS} and \ref{general coprod and comult},
it is straightforward to generalize Lemma
\ref{LS prod geom discription} to the
case when $\lan\sigma,\tau\ran$ doesn't act on $I$ transitively.
\end{rem}

\section{{\bf The wreath product
orbifolds}}\label{sec:wreathproductorbifold}
\begin{defn}
Let $X$ be a compact almost complex manifold with an action $\rho$
of $G$. Let $I$ be a finite set of cardinality $n$. There is a
natural right action of the wreath product $G^I \rtimes \Sigma_I$ on
$X^I$, which we also denote by $\rho$. Namely, $\rho(g\sigma)x \in
X^I$ is defined by
$(\rho(g\sigma)x)_i:=\rho(g_{\sigma(i)})x_{\sigma(i)}$ for all
$g\sigma \in G^I \rtimes \Sigma_I$. Thus, we have an orbifold $[X^I/
\GSI]$ which we call a $\emph{wreath product orbifold}$.
\end{defn}

The following lemma is due to \cite{WZ}.

\begin{lem}\label{prop:locus}
Choose $i_a \in a$ for each $a\in o(\sigma)$. For all $g \in G^I$,
let $\g:=\psi^{\sigma}(g)$ be the cycle product defined in Equation
(\ref{CycleProduct}). The fixed point locus of $g\sigma$ in $X^I$
satisfies
\[
(X^I)^{g\sigma} = \rho(\nu^{\sigma}(g)^{-1}) \prod_{a \in
o(\sigma)}\Delta^a_{X^{\g_a}}
\]
where $\nu^{\sigma}(g) \in G^I$ is defined in Equation (\ref{nu}).
\end{lem}

\Pf It suffices to show the lemma when $\sigma$ is a full
cyclic permutation. In that case, we have $o(\sigma)=\{I\}$. Let
$i$ be the chosen element in $I$ so that
$\g=g_{\sigma^{n-1}(i)}\cdots g_{\sigma(i)}g_i$.
Let $\epsilon_{\g}$ be the element in $G^I$ defined in Equation
(\ref{def-epsilon}). For each $j\in I$, we have
\[
(\rho(\epsilon_{\g}\sigma)x)_j=
\begin{cases}
\rho(\g) x_i & \text{ if $\sigma(j)=i$ } \\
x_{\sigma(j)} & \text{ otherwise. }
\end{cases}
\]
Therefore, $(X^I)^{\epsilon_{\g}\sigma}=\Delta^I_{X^{\g}}$. On the
other hand, we have
$g\sigma=\nu^{\sigma}(g)\epsilon_{\g}\sigma\nu^{\sigma}(g)^{-1}$ by
Equation (\ref{ngn=epsilon}). Hence,
\[
(X^I)^{g\sigma}=\rho(\nu^{\sigma}(g)^{-1})(X^I)^{\epsilon_{\g}\sigma}
=\rho(\nu^{\sigma}(g)^{-1})\Delta^I_{X^{\g}}.
\]
 \qed
\

Choose a representative for each conjugacy class in $\overline{G}$
once and for all. For any index set $J$ and $\overline{\g} \in
\overline{G}^J$, let $\g$ be the representative of $\overline{\g}$
such that $\g_j$ is the chosen representative of the conjugacy class
$\overline{\g}_j$ for each $j \in J$. Let $\cH(X^I, G^I \rtimes
\Sigma_I)$ be the stringy cohomology of  the $(\GSI)$-space $X^I$
introduced in Section \ref{Orbifold}. By Proposition \ref{prop:orb},
the $G^I$-coinvariants of the stringy cohomology is
\[
\cH(X^I, G^I \rtimes \Sigma_I)^{G^I}=\bigoplus_{\sigma \in \Sigma_I}
\bigoplus_{\overline{\g} \in \overline{G}^{o(\sigma)}} \left(
\bigoplus_{g\sigma \in \cO_{\g}\sigma} H^{*} ((X^I)^{g\sigma})
\right)^{G^I}
\]
On the other hand, the Lehn-Sorger algebra associated to
$H_{orb}^{*}([X/G])$ is
\[
H_{orb}^{*}([X/G])\{\Sigma_I\} = \bigoplus_{\sigma \in \Sigma_I}
\bigoplus_{\overline{\g} \in \overline{G}^{o(\sigma)}} \left(
\bigoplus_{\g'\in \overline{\g}} H^{*} ((X^{o(\sigma)})^{\g'})
\right)^{G^{o(\sigma)}} \cdot \sigma.
\]
\begin{prop}\label{linear-iso}
There is a canonical isomorphism of $\Sigma_I$-graded
$\Sigma_I$-modules:
\[
H_{orb}^{*}([X/G])\{\Sigma_I\}  \iso \cH(X^I, G^I \rtimes
\Sigma_I)^{G^I}.
\]
\end{prop}

\Pf Choose $i_a \in a$ for each $a \in o(\sigma)$. On the left-hand
side, we can write
\[
\left( \bigoplus_{\g'\in \overline{\g}} H^{*}
((X^{o(\sigma)})^{\g'})\right)^{G^{o(\sigma)}}\cdot\sigma =\left\{
\left. \left(\frac{1}{|G|^{o(\sigma)}}\sum_{\f\in
G^{o(\sigma)}}\rho(\f)_{*}\bx\right)\sigma \ \right|\ \bx \in
H^{*}((X^{o(\sigma)})^{\g})^{\cZ_{G^{o(\sigma)}}(\g)}
 \right\}
\]
where $\rho(\f): (X^{o(\sigma)})^{\g} \iso
(X^{o(\sigma)})^{\f^{-1}\g\f}$. On the right-hand side,
\[
\left( \bigoplus_{g\sigma \in \cO_{\g\sigma}} H^{*}
((X^I)^{g\sigma})\right)^{G^I} =\left\{\left.
\frac{1}{|G|^{o(\sigma)}} \sum_{f \in G^I} \rho(f)_{*}\bv \ \right|\
\bv \in
H^{*}((X^I)^{\epsilon_{\g}\sigma})^{\cZ_{G^I}(\epsilon_{\g}\sigma)}
\right\}
\]
where $\rho(f):(X^I)^{\epsilon_{\g}\sigma } \iso
(X^I)^{f^{-1}\epsilon_{\g}\sigma f}$. We have
\[
(X^{o(\sigma)})^{\g}= \prod_a X^{\g_a} \cong \prod_a
\Delta^a_{X^{\g_a}} = (X^I)^{\epsilon_{g}\sigma}
\]
where the second isomorphism is defined by the diagonal embedding
which is equivariant with respect to the actions of $\cZ_{G^{o(\sigma)}}(\g)$
and $\cZ_{G^I}(\epsilon_{\g}\sigma)$. For $\bx\in
H^{*}((X^{o(\sigma)})^{\g})^{\cZ_{G^{o(\sigma)}}(\g)}$, let
$\Delta_{\bx}$ denote the corresponding element in
$H^{*}((X^I)^{\epsilon_{\g}\sigma})^{\cZ_{G^I}(\epsilon_{\g}\sigma)}$.
The isomorphism in the proposition is defined by
\[
\frac{1}{|G|^{o(\sigma)}}\sum_{\f \in G^{o(\sigma)}}
\rho(\f)_{*}\bx\sigma \mapsto \frac{1}{|G|^{o(\sigma)}}\sum_{f \in
G^I} \rho(f)_{*}\Delta_{\bx}.
\]
Since we are averaging over $G^{o(\sigma)}$ and $G^I$, this map is
independent of the choice of representatives of conjugacy classes
and of the choice of $\{i_a\}_{a \in o(\sigma)}$. \qed

\section{{\bf The obstruction bundle of the wreath product orbifold}}
\label{sec:obstruction bundle}

In this section, we compute the obstruction bundle $\sR$ introduced
in Section \ref{Orbifold}, for the stringy cohomology of
$(\GSI)$-space $X^I$.

We henceforward adopt the following notation.
\begin{defn}
For all $g\sigma$ and $h\tau$ in $\GSI$, let $\bS_{g\sigma}$ be the
class $\sS_{g\sigma}$ in $K\left((X^I)^{g\sigma}\right)$ defined by
Equation $(\ref{def-Sbundle})$ and let
$\bc(g\sigma,h\tau)$ be the top Chern class
of the obstruction bundle $\sR(g\sigma,h\tau)$.
If $C$ is a connected component of $(X^I)^{g\sigma}$,
let $\mathbf{age}(g\sigma)_C$ denote the age of $g\sigma$ on $C$.
\end{defn}

The following theorem is crucial in proving the algebra isomorphism in
Theorem \ref{ringiso}.

\begin{thm}\label{S-bundle}
Let $\sigma \in \Sigma_I$ and choose $i_a \in a$ for each $a \in
o(\sigma)$. Let $\g \in G^{o(\sigma)}$. Let $\epsilon_{\g}\sigma$ be
the element defined by Equation (\ref{def-epsilon}). We have
\[
\bS_{\epsilon_{\g}\sigma}=\prod_{a \in o(\sigma)} \left(
\Delta^a_{*} \left( \sS_{\g_a} \oplus \frac{|a|-1}{2} TX
|_{X^{\g_a}}\right) \right)
\]
where $\Delta^a : X^{\g_a} \iso \Delta^a_{X^{\g_a}}$ is the
restriction of the diagonal embedding $\Delta^a:X \inc X^a$ and
$\sS_{\g_a}$ is the class in $K(X^{\g_a})$ defined by Equation
(\ref{def-Sbundle}) with respect to the action of $G$ on $X$. For
all $g \in G^I$, we have
\[
\bS_{g\sigma} = \rho(\nu^{\sigma}(g))^{*} \bS_{\epsilon_{\g}
\sigma},
\]
where $\g:=\psi^{\sigma}(g)$.
\end{thm}

\Pf Since $\nu^{\sigma}(g)^{-1}g\sigma\nu^{\sigma}(g)
= \epsilon_{\g}\sigma$ as in Equation (\ref{ngn=epsilon}),
the second claim follows from the first claim and Equation
(\ref{equiv S}).
To prove the first claim, we can assume that
$\sigma$ is a full cyclic permutation without loss of generality. In
that case, $o(\sigma)=\{I\}$  and choose a
representative $j_0 \in I$.
Let $\Delta:X \rightarrow X^I$ be the diagonal map.

Let $V:=\C^I$ be the representation of $\lan\sigma\ran$ induced by
the natural action of $\Sigma_I$ on $\C^I$. Let $\{\be_j\}_{j \in
I}$ be a basis of $\C^I$ such that the action of $\sigma$ on $V$
which we also denote by $\rho$, is
\[
\rho(\sigma)\left(\sum_{j \in I}v_j\be_j\right)=\sum_{j \in
I}v_j\be_{\sigma^{-1}(j)}
\]
for every $\bv=\sum_{j \in I}v_j\be_j \in V$. As a
$\lan\epsilon_{\g}\sigma\ran$-equivariant vector bundle, $TX^I
|_{(X^I)^{\epsilon_{\g} \sigma}}$ is isomorphic to $(T\Delta_X
\otimes V ) |_{\Delta_{X^{\g}}}$. If $p \in \Delta_{X^{\g}}$ and
$\bu\otimes\bv \in T_p \Delta_X \otimes V$, then
$\epsilon_{\g}\sigma$
 acts on $\bu\otimes\bv$ as follows:
\begin{eqnarray*}
\rho(\epsilon_{\g}\sigma)( \bu\otimes\bv )
&=&\rho(\sigma)\left(\rho(\g)\bu \otimes v_i \be_i
+\sum_{j \not= j_0} \bu \otimes v_j \be_j \right)\\
&=&\rho(\g)\bu \otimes v_{j_0}\be_{\sigma^{-1}(j_0)} + \sum_{j \not=
j_0} \bu \otimes v_j \be_{\sigma^{-1}(j)}.
\end{eqnarray*}
Let $r$ be the order of $\g$ and let
$T\Delta_X|_{\Delta_{X^{\g}}}=\bigoplus_{l=0}^{r-1} U_l $ be the
eigenbundle decomposition of the diagonal action of $\g$ where the
eigenvalue of $\rho(\g)$ on the eigenbundle $U_l$ is $\exp{(-2\pi i
l/r)}$. Let $V=\bigoplus_{k=0}^{n-1}V_k$ be the eigenspace
decomposition of $\sigma$ on $V$ where $n$ is the cardinality of
$I$. The eigenvalue of $\rho(\sigma)$ on $V_k$ is $\exp{(-2\pi
ik/n)}$. If $V_k$ is generated by $\bv_k=\sum_{j \in I} v_{k,j}
\be_j$, then the equality $\rho(\sigma)\bv_k=\exp(-2\pi ik/n)\bv_k$
implies
\[
\sum_{j \in I} v_{k,j} \be_{\sigma^{-1}(j)}=\sum_{j \in I}\exp(-2\pi
ik/n)v_{k,j} \be_j.
\]
By comparing the coefficient of $\be_j$, we obtain
\begin{equation}\label{vsj=evj}
v_{k, \sigma(j)}=\exp(-2\pi i k/n) v_{k,j}.
\end{equation}
Let $V^l_k$ is a 1-dimensional subspace spanned by
\[
 \bv^l_k = \sum_{m=0}^{n-1} \exp(-2\pi i
\frac{l m}{n r})v_{k,\sigma^m(j_0)} \be_{\sigma^m(j_0)}.
\]
Introduce another decomposition $V=\bigoplus_{k=0}^{n-1}V^l_k$, and
then, together with the decomposition
$T\Delta_X|_{\Delta_{X^{\g}}}=\bigoplus_{l=0}^{r-1} U_l$, we have
\begin{equation}\label{eq:decom}
TX^I|_{(X^I)^{\epsilon_{\g}\sigma}}=\bigoplus_{k=0}^{n-1} \left(
\bigoplus_{l=0}^{r-1}  U_l \otimes V^l_k  \right).
\end{equation}
This turns out to be the eigenbundle decomposition of the action of
$\epsilon_{\g}\sigma$ on $TX^I |_{(X^I)^{\epsilon_{\g} \sigma}}$
and the eigenvalue of $U_l \otimes V^l_k$ is $\exp(-2\pi
i\left(\frac{l}{nr}+ \frac{k}{n}\right))$. In fact, for any $\bu_l
\in U_l$,
\begin{eqnarray*}
\rho(\epsilon_{\g}\sigma)_{*} \bu_l\otimes\bv^l_k
&=&\rho(\sigma)_{*}\left(\rho(\g)_{*}\bu_l\otimes v_{k,j_0}
\be_{j_0} + \sum_{m=1}^{n-1} \bu_l\otimes \exp\left(-2\pi i \frac{l m
}{n r}\right)v_{k,\sigma^m(j_0)} \be_{\sigma^m(j_0)}
 \right) \\
&=&\exp\left(-2\pi i \frac{l}{r}\right)\bu_l \otimes v_{k,j_0}
\be_{\sigma^{-1}(j_0)} +\sum_{m=1}^{n-1}\bu_l\otimes \exp\left(-2\pi
i \frac{l m }{n
r}\right)v_{k,\sigma^m(j_0)}\be_{\sigma^{m-1}(j_0)} \\
&=&\exp\left(-2\pi i \frac{l}{r}\right)\bu_l \otimes \exp\left(-2\pi i
\frac{k}{n}\right)v_{k,\sigma^{-1}(j_0)} \be_{\sigma^{-1}(j_0)}\\
&& \ \ \ \ +\sum_{m=1}^{n-1} \bu_l \otimes \exp\left(-2\pi i
\left(\frac{l m}{nr} +\frac{k}{n}\right)\right)v_{k,\sigma^{m-1}(j_0)}
\be_{\sigma^{m-1}(j_0)}\\
&=&\bu_l \otimes \left(\sum_{m=0}^{n-1}  \exp\left(-2\pi i
\left(\frac{l m+l}{nr}
+\frac{k}{n}\right)\right)v_{k,\sigma^{m}(j_0)}\be_{\sigma^{m}(j_0)}
 \right)\\
&=& \exp\left(-2\pi i \left(\frac{l}{nr}+ \frac{k}{n}\right)\right)
\bu_l \otimes \bv^l_k
\end{eqnarray*}
where the first and second equalities follow from the definition of
the action of $\rho(\epsilon_{\g})$ and $\rho(\sigma)$, and the
third equality follows from Equation (\ref{vsj=evj}). Thus we have
\[
\bS_{\epsilon_{\g}\sigma}=\bigoplus_{k=0}^{n-1}
\bigoplus_{l=0}^{r-1} \left(\frac{l}{nr}+ \frac{k}{n}\right) U_l
\otimes V^l_k  \ .
\]
For all $k$ and $l$, we have $U_l\otimes V^l_k\cong U_l$ since $V^l_k$ is
a 1-dimensional vector space. Thus
\[
\bS_{\epsilon_{\g}\sigma}=\bigoplus_{k=0}^{n-1}
\bigoplus_{l=0}^{r-1} \left(\frac{l}{nr}+ \frac{k}{n}\right) U_l =
\bigoplus_{l=0}^{r-1} \frac{l}{r}U_l \oplus \bigoplus_{k=0}^{n-1}
\frac{k}{n} T\Delta_X |_{\Delta_{X^{\g_a}}} =
\Delta_{*}\left(\sS_{\g} \oplus \frac{n-1}{2}TX|_{X^{\g}}\right).
\]
\qed

This theorem leads to the following corollary which was obtained
in \cite{WZ} through the direct calculation.

\begin{cor}\label{age}
Let $\g:=\psi^{\sigma}(g)$ and let $C_{\g_a}$ be a connected
component of $X^{\g_a}$ for each $a\in o(\sigma)$. Every connected
component of $(X^I)^{g\sigma}$ can be written as
$C:=\rho(\nu^{\sigma}(g))^{-1}\left(\prod_a
\Delta^a_{C_{\g_a}}\right)$ and we have
\[
\mathbf{age}(g\sigma)_C = \frac{\dim X \cdot l_{\sigma}}{2}+ \sum_{a
\in o(\sigma)} age(\g_a)_{C_{\g_a}},
\]
where $age(\g_a)_{C_{\g_a}}$ is the age of $\g_a$ on ${C_{\g_a}}$
with respect to the action $G$ on $X$ and $l_{\sigma}$ is the length
of $\sigma$.
\end{cor}

\emph{For the rest of the section, we assume that $G$ is Abelian and
 adopt the following notation.}
\begin{notation}
Let $\sigma,\tau \in \Sigma_I$. Let $a \in o(\sigma), b \in o(\tau),
c \in o(\sigma\tau)$ and $d \in o(\sigma,\tau)$. Let
\begin{eqnarray*}
o(\sigma)_d&:=&\{a \in o(\sigma) \ |\ a \subset d  \},\\
o(\tau)_d&:=&\{b \in o(\tau) \ |\ b \subset d  \}, \\
o(\sigma\tau)_d&:=&\{c \in o(\sigma\tau) \ |\ c \subset d  \}.
\end{eqnarray*}
Once and for all, choose representatives $i_a \in a$, $i_b \in b$,
and $i_c \in c$ for all $a \in o(\sigma)$, $b \in o(\tau)$, and $c
\in o(\sigma\tau)$. Furthermore, let $gd(d):=gd(\sigma,\tau)_d$ for
all $d \in o(\sigma,\tau)$. Let $f_d$ be the image of $f \in G^I$ by
the obvious projection from
$G^I$ to $G^d$. If $\g\in G^{o(\sigma)}$, let $\g_d$ be the image of $\g$ by
the obvious projection from $G^{o(\sigma)}$ to $G^{o(\sigma)_d}$.
Define $\h_d, \w_d$ in the same manner for all
$\h \in G^{o(\tau)}$ and $\w \in G^{o(\sigma\tau)}$.
\end{notation}

\begin{lem}\label{intersection}
For all $\g \in G^{o(\sigma)}, \h\in G^{o(\tau)}, \w \in
G^{o(\sigma\tau)}$ and $f \in G^I$ such that
$\epsilon_{\g}\sigma\cdot f^{-1}\epsilon_{\h}\tau f$ lies in
$\cO_{\w}\sigma\tau$, there exists $\f_{(d)} \in G^{2gd(d)}$ for each $d
\in o(\sigma,\tau)$ and $\fbar \in \prod_{a} \Delta^a_G$ such that
\begin{equation}\label{intersection eq}
 \prod_{a \in o(\sigma)} \Delta^a_{X^{\g_a}} \cap
\rho(f)\prod_{b \in o(\tau)}\Delta^b_{X^{\h_b}} = \prod_{d \in
o(\sigma,\tau)} \rho(\fbar_d)^{-1}\Delta_{X^{\g_d, \h_d, \w_d,
\f_{(d)}}}.
\end{equation}
\end{lem}
Note that $(X^I)^{\epsilon_{\g}\sigma}=\prod_{a \in o(\sigma)} \Delta^a_{X^{\g_a}}$
and $(X^I)^{f^{-1}\epsilon_{\h}\tau f}=\rho(f)\prod_{b \in o(\tau)}\Delta^b_{X^{\h_b}}$
by Lemma \ref{prop:locus}.

\Pf  Since the left-hand side of Equation (\ref{intersection eq})
breaks up into the direct product
\[
\prod_d \left( \prod_{a \in o(\sigma)_d} \Delta^a_{X^{\g_a}} \cap
\rho(f_d)\prod_{b \in o(\tau)_d}\Delta^b_{X^{\h_b}}\right),
\]
we can assume that $\lan\sigma,\tau\ran$ acts transitively on $I$
without loss of generality. Let $gd:=gd(I)$ and $\f:=\f_{(I)} \in G^{2gd}$.
Since $\lan\sigma,\tau\ran$ acts on $I$ transitively, the
intersection $\prod_a\Delta^a_{X^{\g_a}} \cap \rho(f)
\prod_b\Delta^b_{X^{\h_b}}$ is contained in
$\rho(\fbar)^{-1}\Delta_{X^{\g,\h}}$ for some $\fbar \in \prod_a\Delta^a_G$.

Associate an unoriented graph $\Gamma$ to $\sigma$ and $\tau$, where
the vertices of $\Gamma$ are the elements of $I$ and the edges of
$\Gamma$ are $\{\sigma^{k_a}(i_a), \sigma^{k_a+1}(i_a)\}$ and
$\{\tau^{k_b}(i_b), \tau^{k_b+1}(i_b)\}$ for all $a \in o(\sigma)$,
$b \in o(\tau)$, $k_a=0,\cdots,|a|-2$, and $k_b=0,\cdots,|b|-2$.
This graph $\Gamma$ is connected since $\lan\sigma,\tau\ran$ acts on
$I$ transitively. The Euler characteristic of $\Gamma$ is
$n-l_{\sigma}-l_{\tau}$. If $b_1$ is the first Betti number of
$\Gamma$, then $1-b_1=n-l_{\sigma}-l_{\tau}$, hence
\begin{equation}\label{betti}
b_1=l_{\sigma}+l_{\tau}+1-n=2gd+|o(\sigma\tau)|-1.
\end{equation}
Take $z \in \rho(\fbar)^{-1}\Delta_{X^{\g,\h}}$
and then $z$ is in the intersection $\prod_a\Delta^a_{X^{\g_a}} \cap
\rho(f)\prod_b\Delta^b_{X^{\h_b}}$ if and only if $z$
satisfies, for every edge $\{v_0,v_1\}$,
\[
z_{v_0}=
\begin{cases}
z_{v_1}&\text{ if
$\{v_0,v_1\}=\{\sigma^{k_a}(i_a),\sigma^{k_a+1}(i_a)\}$ for some $a$
and $k_a$,}\\
\rho(f_{v_0}f_{v_1}^{-1})z_{v_1} & \text{ if
$\{v_0,v_1\}=\{\tau^{k_b}(i_b), \tau^{k_b+1}(i_b)\}$ for some $b$
and $k_b$.}
\end{cases}
\]

Let $\alpha$ be a closed, oriented circle in the graph and let
$E_{\alpha}$ be the set of oriented edges of $\Gamma$ which are
contained in $\alpha$ and whose orientations are induced from the
orientation of $\alpha$. The oriented edge associated to the edge
$\{v_0,v_1\}$ is denoted by $(v_0,v_1)$. Let
\[
f_{v_0,v_1}:=
\begin{cases}
1&\text{ if
$\{v_0,v_1\}=\{\sigma^{k_a}(i_a),\sigma^{k_a+1}(i_a)\}$ for some $a$
and $k_a$,}\\
f_{v_0}f_{v_1}^{-1}& \text{ if
$\{v_0,v_1\}=\{\tau^{k_b}(i_b), \tau^{k_b+1}(i_b)\}$ for some $b$
and $k_b$.}
\end{cases}
\]
and let $\f_{\alpha}:=\prod_{(v_0,v_1) \in E_{\alpha}}f_{v_0,v_1}$,
then $z_i=\rho(\f_{\alpha})z_i$ for every vertex $i$ in $\alpha$.
Hence, $z$ is in the intersection $\Delta^a_{X^{\g_a}} \cap
\rho(f)\prod_{b \in o(\tau)}\Delta^b_{X^{\h_b}}$ if and only if $z
\in \rho(\fbar)^{-1}\Delta_{X^{\g,\h, \f_{\alpha}}}$ for every circle
$\alpha$. If $\overline{\alpha}$ denotes the same circle $\alpha$
with the opposite orientation, then
$\f_{\overline{\alpha}}=\f_{\alpha}^{-1}$. If $\alpha$ is homologous
to $\alpha'$, then $\f_{\alpha'}=\f_{\alpha}$. Therefore, if
$\{\alpha_k\}_{k=1,\cdots,b_1}$ is a basis of $H_1(\Gamma, \Z)$, we
have
\[
\prod_{a \in o(\sigma)} \Delta^a_{X^{\g_a}} \cap \rho(f)\prod_{b \in
o(\tau)}\Delta^b_{X^{\h_b}}=\rho(\fbar)^{-1}\Delta_{X^{\g, \h, \f'}}
\] where $\f'
\in G^{2gd+|o(\sigma\tau)|-1}$ and $\f'_k:=\f_{\alpha_k}$.
Furthermore, since $X^{\g,\h,\f'}$ has to be contained in
$X^{\w}$, we have $X^{\g,\h,\f'}=X^{\g,\h,\w,\f}$ for some $\f\in
G^{2gd}$. \qed

\begin{prop}\label{R(gs,ht)}
Let $\g \in G^{o(\sigma)}$, $\h \in G^{o(\tau)}$, $\w \in
G^{o(\sigma\tau)}$ and $f\in G^I$ such that
$\epsilon_{\g}\sigma\cdot f^{-1}\epsilon_{\h} \tau f \in
\cO_{\w}\sigma\tau$. Let $Z_d:= X^{\g_d,\h_d,\w_d}$. Choose $\f_{(d)}
\in G^{2{gd(d)}}$ for each $d \in o(\sigma,\tau)$ and $\fbar  \in
\prod_a \Delta^a_G$ which satisfies the equality (\ref{intersection
eq}) in Lemma \ref{intersection}. We have
\begin{equation}\label{R(egs,eht)}
\sR(\epsilon_{\g}\sigma, f^{-1}\epsilon_{\h} \tau f)= \prod_{d \in
o(\sigma,\tau)} \rho(\fbar_d )^{*}\Delta^d_* \left( TZ_d^{\f_{(d)}}
\oplus (gd(d)-1) TX|_{Z_d^{\f_{(d)}}} \oplus \sS_{\g_d}|_{Z_d^{\f_{(d)}}}
\oplus\sS_{\h_d}|_{Z_d^{\f_{(d)}}} \oplus\sS_{\w_d^{-1}}|_{Z_d^{\f_{(d)}}}
\right)
\end{equation}
where $\Delta^d:Z^{\f_{(d)}} \iso \Delta^d_{Z_d^{\f_{(d)}}}$ is the
isomorphism induced by the diagonal embedding $X \inc X^d$.

For any $g,h\in G^I$,
\[
\sR(g\sigma,h\tau)
=\rho(\nu^{\sigma}(g))^{*}\sR(\epsilon_{\g}\sigma,f^{-1}\epsilon_{\h}\tau{f})
\]
where $\g:=\psi^{\sigma}(g)$, $\h:=\psi^{\tau}(h)$, and
$f:=\nu^{\tau}(h)^{-1}\nu^{\sigma}(g)$.
\end{prop}

\Pf The second claim simply follows from Equation
(\ref{ngn=epsilon}) and (\ref{G-equivR}). Let $\sigma_d$ and
$\tau_d$ be the restriction of the action of $\sigma$ and $\tau$ on
$X^d$ respectively. The left-hand side of Equation
(\ref{R(egs,eht)}) breaks up into the external direct product
\[
\prod_d \sR_d((\epsilon_{\g})_d\sigma_d, f_d^{-1}(\epsilon_{\h})_d
\tau_d f_d)
\]
where $\sR_d$ denotes the obstruction bundle of $(G^d\rtimes
\Sigma_d)$-space $X^d$. Hence, we can assume that
$\lan\sigma,\tau\ran$ acts transitively on $I$ without loss of
generality. Let $gd:=gd(I)$.

Let $f' \in G^I$ such that
\[
\epsilon_{\g}\sigma\cdot f^{-1}\epsilon_{\h} \tau f
\cdot f'^{-1}\epsilon_{\w^{-1}}(\sigma\tau)^{-1}f'=1
\]
By Equation (\ref{G-equivR}),
\[
\rho(\fbar)_{*}\sR(\epsilon_{\g}\sigma, f^{-1}\epsilon_{\h} \tau f)
=\sR(\fbar^{-1}\epsilon_{\g}\sigma\fbar, (f \fbar)^{-1}\epsilon_{\h}
\tau f\fbar)
\]
which is equal to
\[
T(\Delta_{Z^{\f}}) \ominus TX^I |_{\Delta_{Z^{\f}}}
\oplus\bS_{\fbar^{-1}\epsilon_{\g}\sigma\fbar}|_{\Delta_{Z^{\f}}}\oplus
\bS_{(f\fbar)^{-1}\epsilon_{\h} \tau f\fbar} |_{\Delta_{Z^{\f}}}
\oplus
\bS_{(f'\fbar)^{-1}\epsilon_{\w^{-1}}(\sigma\tau)^{-1}f'\fbar}
|_{\Delta_{Z^{\f}}}.
\]
Since $\fbar$ commutes with $\epsilon_{\g}\sigma$, we have
\[
\bS_{\fbar^{-1}\epsilon_{\g}\sigma\fbar} =\bS_{\epsilon_{\g}\sigma}.
\]
Since the commutator $[(f\fbar)^{-1},\epsilon_{\h} \tau]$ belongs to
$\lan\g,\h,\w,\f\ran$, the actions of $(f\fbar)^{-1}\epsilon_{\h}
\tau f\fbar$ and $\epsilon_{\h}\tau$ coincide on $\Delta_{Z^{\f}}$.
Therefore
\[
\bS_{(f\fbar)^{-1}\epsilon_{\h}\tau f\fbar}|_{\Delta_{Z^{\f}}}
=\bS_{\epsilon_{\h} \tau } |_{\Delta_{Z^{\f}}} .
\]
Since also the commutator $[(f'\fbar)^{-1},\epsilon_{\w^{-1}}
(\sigma\tau)^{-1}]$ belongs to $\lan\g,\h,\w,\f\ran$, the actions of
$\epsilon_{\w^{-1}}(\sigma\tau)^{-1}$ and
$(f\fbar)^{-1}\epsilon_{\w^{-1}} (\sigma\tau)^{-1} f\fbar$ coincide
on $\Delta_{Z^{\f}}$. Therefore
\[
\bS_{(f'\fbar)^{-1}\epsilon_{\w^{-1}}(\sigma\tau)^{-1}f'\fbar}
|_{\Delta_{Z^{\f}}} =\bS_{\epsilon_{\w^{-1}}(\sigma\tau)^{-1}}
|_{\Delta_{Z^{\f}}}.
\]
Thus
\[\rho(\fbar)_{*}\sR(\epsilon_{\g}\sigma, f^{-1}\epsilon_{\h} \tau f)=
T(\Delta_{Z^{\f}}) \ominus TX^I |_{\Delta_{Z^{\f}}} \oplus
\bS_{\epsilon_{\g}     \sigma           }|_{\Delta_{Z^{\f}}} \oplus
\bS_{\epsilon_{\h}     \tau             }|_{\Delta_{Z^{\f}}} \oplus
\bS_{\epsilon_{\w^{-1}}(\sigma\tau)^{-1}}|_{\Delta_{Z^{\f}}}
\]
and the proposition follows from Theorem \ref{S-bundle}. \qed

\begin{defn}\label{fc[ghwf]}
Let $\fc[\g,\h,\w,\f]_d$ be the top Chern class of the bundle
\[
\left( TZ_d^{\f_{(d)}} \oplus (gd(d)-1) TX|_{Z_d^{\f_{(d)}}} \oplus
\sS_{\g_d}|_{Z_d^{\f_{(d)}}} \oplus\sS_{\h_d}|_{Z_d^{\f_{(d)}}}
\oplus\sS_{\w_d^{-1}}|_{Z_d^{\f_{(d)}}}\right)
\]
and let $\fc[\g,\h,\w,\f]:=\bigotimes_d \fc[\g,\h,\w,\f]_d$.
\end{defn}
\begin{cor}
The top Chern class of $\sR(\epsilon_{\g}\sigma,
f^{-1}\epsilon_{\h}\tau f)$ is
\[
\bc(\epsilon_{\g}\sigma, f^{-1}\epsilon_{\h}\tau
f)=\bigotimes_{d}\rho(\fbar_d )^{*}\Delta^d_{*}\fc[\g,\h,\w,\f]_d
=\rho(\fbar )^{*}\Delta_{*} \fc[\g,\h,\w,\f]
\]
where $\Delta=\prod_d \Delta_d :\prod_d Z^{\f_{(d)}} \iso \prod_d
\Delta^d_{Z_d^{\f_{(d)}}}$.
\end{cor}

\section{{\bf The ring isomorphism}}\label{sec:ring isomorphism}

Suppose that $G$ is an Abelian group and that $\lan\sigma,\tau\ran$
acts transitively on $I$. Let $gd:=gd(\sigma,\tau)_I$. Let $\g \in
G^{o(\sigma)}$, $\h \in G^{o(\tau)}$, $\w \in G^{o(\sigma\tau)}$ and
$f\in G^I$ such that $\epsilon_{\g}\sigma\cdot f^{-1}\epsilon_{\h}
\tau f \in \cO_{\w}\sigma\tau$. Let $Z:= X^{\g,\h,\w}$. Choose
$\f\in G^{2{gd}}$ and $\fbar\in\prod_a \Delta^a_G$ so that
\[\prod_{a \in o(\sigma)}\Delta^a_{X^{\g_a}}\cap\rho(f)\prod_{b\in
o(\tau)}
\Delta^b_{X^{\h_b}}=\rho(\fbar)^{-1}\Delta_{X^{\g,\h,\w,\f}}\] as in
Lemma \ref{intersection}. Let $\bp_{\w},\bq_{\w}$ be the following
diagonal embeddings
\[
\prod_a X^{\g_a} \times \prod_b X^{\h_b}
\stackrel{\bp_{\w}}{\longleftarrow} Z
\stackrel{\bq_{\w}}{\longrightarrow} \prod_c X^{\w_c}
\]
and $\br_{\f}:Z^{\f} \inc Z$ be the canonical inclusion. For all
$\bx \in H^{*}\left(\prod_a X^{\g_a}\right)^{G^{o(\sigma)}}$, let
$\Delta_{\bx} \in
H^{*}\left((X^I)^{\epsilon_{\g}\sigma}\right)^{\prod_a \Delta^a_G}$
be the push-forward of $\bx$ by the isomorphism $\prod_a X^{\g_a}
\iso \prod_a \Delta_{X^{\g_a}}$.

\begin{lem}\label{MainLemma}
Under the canonical isomorphism in Proposition \ref{linear-iso},
the element
\begin{equation}\label{7-1-a}
\frac{1}{|G|^{|o(\sigma\tau)|}}\sum_{f' \in G^I}
\rho(f')_*\left(\Delta_{\bx} \cdot \rho(f)_*\Delta_{\by}\right) \in
\left(\bigoplus_{w\sigma\tau \in \cO_{\w}\sigma\tau} H^{*}
((X^I)^{w\sigma\tau})\right)^{G^I}
\end{equation}
 corresponds to
\begin{equation}\label{7-1-b}
\bq_{\w * }\left(\bp^{*}_{\w} (\bx  \otimes \by )\cup
\br_{\f*}\fc[\g,\h,\w,\f]\right)\cdot \sigma\tau \in H^{*} (
(X^{o(\sigma\tau)})^{\w})^{G^{o(\sigma\tau)}}\sigma\tau.
\end{equation}
\end{lem}
\Pf By the change of variables $f'=\fbar f''$,
\[
\frac{1}{|G|^{|o(\sigma\tau)|}}\sum_{f' \in G^I}
\rho(f')_*(\Delta_{\bx} \cdot \rho(f)\Delta_{\by})=
\frac{1}{|G|^{|o(\sigma\tau)|}}\sum_{f' \in G^I}
\rho(f'')_*\left(\rho(\fbar)_*\Delta_{\bx}\cdot
\rho(f\fbar)_*\Delta_{\by} \right).
\]
Since $\Delta_{\bx}$ is $\prod_a \Delta^a_G$-invariant and $\fbar
\in \prod_a \Delta^a_G$, we have
$\rho(\fbar)_*\Delta_{\bx}=\Delta_{\bx}$.
The action of $f\fbar$
restricted to $\Delta_{Z^{\f}}$ agrees with the action
of some element $\gamma$ in $\prod_b \Delta^b_G$, since
$\rho(f\fbar)\prod_{b\in o(\tau)}\Delta^b_{X^{\h_b}}$ contains
$\Delta_{X^{\g,\h,\w,\f}}$. This implies that we have
the following commutative diagrams,
\[
\begin{CD}
\Delta_{Z^{\f}} @<\rho(\gamma)<< \rho(\gamma^{-1})\Delta_{Z^{\f}}
@<\rho(\gamma^{-1})<< \Delta_{Z^{\f}} \\
@V \iota_1 VV  @V \iota_2 VV @V \iota_3 VV \\
\rho(\gamma)\prod_b\Delta^b_{X^{\h_b}} @>\rho(\gamma^{-1})>> \prod_b\Delta^b_{X^{\h_b}}
@<\rho((f\fbar)^{-1})<< \rho(f\fbar)\prod_b\Delta^b_{X^{\h_b}}\\
\end{CD}
\]
where $\iota_1, \iota_2, \iota_3$ are the canonical inclusions.
By the diagrams above, we have
\[
\rho(f\fbar)_*(\Delta_{\by})|_{\Delta_{Z^{\f}}}=
\iota_3^*\circ\rho((f\fbar)^{-1})^*\Delta_{\by}=
\iota_1^*\circ\rho(\alpha^{-1})^*\Delta_{\by}=\Delta_{\by}|_{\Delta_{Z^{\f}}},
\]
where the last equality holds because
$\Delta_{\by}$ is $\prod_b \Delta^b_G$-invariant.
Thus, by Equation (\ref{product}) and Proposition
\ref{R(gs,ht)}, Equation (\ref{7-1-a}) is equal to
\begin{equation}\label{77}
\frac{1}{|G|^{|o(\sigma\tau)|}}\sum_{f' \in G^I} \rho(f')_* \left[
\bq_{\Delta_{\w}*} \left( \Delta_{\bx}|_{\Delta_{Z^{\f}}} \cup
\Delta_{\by}|_{\Delta_{Z^{\f}}}\cup
\rho(\fbar)_*\bc(\epsilon_{\g}\sigma, f^{-1}\epsilon_{\h} \tau f)
\right) \right]
\end{equation}
where $\bq_{\Delta_{\w}}: \Delta_{Z^{\f}}\inc\prod_c
\Delta^c_{X^{\w_c}}$ is the diagonal embedding. Therefore, by Definition
\ref{fc[ghwf]}, Equation (\ref{7-1-a}) corresponds to
\[
\bq_{\w*} \circ \br_{\f*} \left( \br_{\f}^{*}\circ \bp^{*}_{\w} (\bx
\otimes \by ) \cup\fc[\g,\h,\w,\f]\right)\sigma\tau
\]
under the isomorphism. We obtain the proposition by applying the
projection formula. \qed

\begin{thm} \label{ringiso}
Assume that $G$ is an Abelian group. $\cH(X^I,\GSI)^{G^I}$ is
canonically isomorphic as a $\Sigma_I$-Frobenius algebra to
$H^{*}_{orb}([X/G])\{\Sigma_I\}$ under the isomorphism defined in
Proposition \ref{linear-iso}. In particular, the Lehn-Sorger algebra
$H^{*}_{orb}([X/G])\{\Sigma_I\}$ satisfies the trace axiom.
\end{thm}

\Pf If we prove that they are isomorphic to each other as rings under
the isomorphism in Proposition \ref{linear-iso}, then all other
properties of $\Sigma_I$-Frobenius algebras are clearly preserved by
the isomorphism and, in particular, the trace axiom on the
Lehn-Sorger side of the equality is satisfied.

Let $\lambda$ be a partition of $I$. Consider the following subspace
of $\cH(X^I, \GSI)$:
\[
\cH(\lambda):=\bigoplus_{o(\sigma) <\lambda} \cH_{\sigma}.
\]
It is clear that this is a subalgebra of $\cH(X^I, \GSI)$. By Lemma
\ref{intersection}, Proposition \ref{R(gs,ht)}, and the K\"{u}nneth
theorem, it is clear that
\[
\cH(\lambda)\cong\bigotimes_{d \in \lambda}\cH(X^d,
G^d\rtimes\Sigma_d).
\]
Since the obstruction bundle is $G^I$-equivariant and the equality
(\ref{intersection eq}) is preserved by the action of $G^I$, we
obtain
\begin{equation}\label{splitting WP}
\cH(\lambda)^{G^I}\cong\bigotimes_{d \in \lambda}\cH(X^d,
G^d\rtimes\Sigma_d)^{G^d}.
\end{equation}
Hence, to prove the theorem, comparing Equation (\ref{splitting LS
product}) with Equation (\ref{splitting WP}), we can assume that
$\lan\sigma, \tau\ran$ acts transitively on $I$ without loss of
generality.

The product of $\cH(X^I,G^I\rtimes \Sigma_I)^{G^I}$ corresponding to
the Lehn-Sorger product $\bx\sigma \cdot \by \tau$ under the
isomorphism is, by the change of variables $f''
=f f'$,
\begin{eqnarray*}
&&\left(\frac{1}{|G|^{|o(\sigma)|}}\sum_{f' \in
G^I}\rho(f')_*\Delta_{\bx}\right)\cdot
\left(\frac{1}{|G|^{|o(\tau)|}}\sum_{f''
\in G^I}\rho(f'')_*\Delta_{\by}\right)\\
&=&\frac{1}{|G|^{|o(\sigma)|+|o(\tau)|-|o(\sigma\tau)|}}\sum_{f \in
G^I}\frac{1}{|G|^{|o(\sigma\tau)|}}
\sum_{f' \in G^I} \rho(f')_* \left(\Delta_{\bx} \cdot \rho(f)_*\Delta_{\by} \right)
\end{eqnarray*}
By Lemma \ref{MainLemma},
$\sum_{f' \in G^I} \rho(f')_* \left(\Delta_{\bx} \cdot \rho(f)_*\Delta_{\by} \right)$
only depends on $\w$ and $\f$. Hence, by the construction of $\f$ in Lemma \ref{intersection}
and by Lemma \ref{surj},
\begin{eqnarray*}
&&\frac{1}{|G|^{|o(\sigma)|+|o(\tau)|-|o(\sigma\tau)|}}\sum_{f \in
G^I}\frac{1}{|G|^{|o(\sigma\tau)|}}
\sum_{f' \in G^I} \rho(f')_* \left(\Delta_{\bx} \cdot \rho(f)_*\Delta_{\by} \right)\\
&=& |G|^{|o(\sigma\tau)|-1}\sum_{\w}\sum_{\f \in G^{2{gd}}}
\frac{1}{|G|^{|o(\sigma\tau)|}}  \sum_{f' \in G^I} \rho(f')_* \left(
\Delta_{\bx} \cdot \rho(f)_*\Delta_{\by} \right),
\end{eqnarray*}
where we chose an $f$ for each pair $(\w,\f)$.
By Lemma \ref{MainLemma}, this corresponds to
\begin{equation}\label{71}
\frac{1}{|G|} \sum_{\fa \in
G^{o(\sigma\tau)}}\rho(\fa)\sum_{\w}\sum_{\f \in G^{2{gd}}}\bq_{ \w
* }  \left(\bp^{*}_{\w}(\bx  \otimes \by ) \cup \br_{\f*}\fc[\g,\h,\w,
\f]\right)\sigma\tau.
\end{equation}
To finish the proof, we need to show Equation $(\ref{71})$ is equal
to Equation $(\ref{61})$,
\[
\bx \sigma\cdot \by \tau =\frac{1}{|G|} \sum_{\fa\in
G^{o(\sigma\tau)}} \rho(\fa) \left[ \sum_{\w} \bq_{\w*} \left(
\bp_{\w}^{*}(\bx \otimes \by) \cup c(\g,\h)|_{Z_{\w}} \cup
\e^{{gd}}|_{Z_{\w}} \cup c(\w^{-1})|_{Z_{\w}} \cup E_{\w} \right)
\right]\cdot \sigma\tau.
\]
By the linearity of $\bq_{\w*}$ and the
cup product, we must show
\begin{equation}\label{final}
\sum_{\f} \br_{\f*}\fc[\g,\h,\w,\f] = c(\g,\h)|_Z \cup c(\w^{-1})
|_Z \cup E_{\w} \cup \e^{{gd}}|_Z.
\end{equation}
The right-hand side of Equation (\ref{final}) is given by
\begin{equation}
c_{top}(\cE) \cup \e^{{gd}}|_Z=
\begin{cases}
c_{top}(\cE) & \text{ if } {gd}=0, \\
c_{top}(\cE) \cup \e & \text{ if } {gd}=1 \text{ and } Z=X, \\
0 & \text{ otherwise }
\end{cases}
\end{equation}
where $\cE$ is the bundle whose class in $K(Z)$ is equal to
\begin{equation}\label{cE}
TZ \ominus TX|_Z \oplus \sS_{\g}|_Z \oplus \sS_{\h}|_Z \oplus
\sS_{\w^{-1}}|_Z.
\end{equation}
Since $\fc[\g,\h,\w,\f]$ is the top Chern class of the bundle
\[
TZ^{\f} \oplus (gd-1) TX|_{Z^{\f}} \oplus \sS_{\g}|_{Z^{\f}}
\oplus\sS_{\h}|_{Z^{\f}} \oplus\sS_{\w^{-1}}|_{Z^{\f}},
\]
the left-hand side of Equation (\ref{final}) is equal to
\[
 \sum_{\f \in G^{2{gd}}} \br_{\f*}
\left[ \br_{\f}^{*} c_{top}(\cE) \cup c_{top} \left( gd\cdot
TX|_{Z^{\f}}\ominus TZ|_{Z^{\f}} \oplus TZ^{\f}\right) \right].
\]
Hence, by the projection formula,
\begin{equation}\label{72}
\sum_{\f} \br_{\f*}\fc[\g,\h,\w,\f] = c_{top}(\cE) \cup \sum_{\f \in
G^{2{gd}}} \br_{\f*}c_{top}\left( gd\cdot TX|_{Z^{\f}}\ominus
TZ|_{Z^{\f}} \oplus TZ^{\f}\right).
\end{equation}
When ${gd}=0$, we have$\sum_{\f} \br_{\f*}\fc[\g,\h,\w,\f]
=c_{top}(\cE)$ since $Z^{\f}=Z$. When ${gd}=1$,
\begin{equation}\label{74}
c_{top}\left(TX|_{Z^{\f}}\ominus TZ|_{Z^{\f}} \oplus TZ^{\f}\right)
= c_{top}(N_{Z/X})|_{Z^{\f_1,\f_2}} \cup c_{top}(TZ^{\f_1,\f_2})
\end{equation}
where $N_{Z/X}$ is the normal bundle of $Z$ in $X$. Hence, the
right-hand side of Equation (\ref{72}) vanishes unless $Z=X$ by
dimensional considerations. When ${gd}=1$ and $Z=X$, we have $\cE=0$ and
the right-hand side of
Equation (\ref{72}) equals $\sum_{\f_1,\f_2}
\br_{\f*} c_{top}\left(TX^{\f_1,\f_2}\right)=\e$. When ${gd} \geq
2$, the right-hand side of Equation (\ref{72}) vanishes by
dimensional considerations. \qed

\section{{\bf Example : a torus with an involution}}\label{sec:ex1}
Let $T:=\C^2/(\Z p_1 +\Z p_2 +\Z p_3 +\Z p_4)$ be a $2$ dimensional
complex torus and $\Z_2$ be a group of order $2$ generated by
$\alpha$. Let $X:=T$ and $G:=\Z_2$ and let $G$ act on $X$ by $\alpha
: (z,w) \mapsto (-z,-w)$. We have an orbifold $[X/G]=[T/\Z_2]$.
There are 16 points, $\{q_j\}_{j=1,\cdots,16}$, in $X^{\alpha}$
corresponding to $\{\Z\cdot \frac{1}{2} p_1 +\Z \cdot \frac{1}{2}p_2
+\Z \cdot \frac{1}{2}p_3 +\Z \cdot \frac{1}{2}p_4\} \subset \C^2$.
The orbifold cohomology of $[X/G]$ is
\[
H^{*}_{orb}([X/G])=H^{*}(T)^{\Z_2}\oplus
H^0(\{q_j\}_{j=1,\cdots,16}).
\]
Let $\phi_1:=\1$ and $\phi_2$ be the class of top dimension such
that $\int_T\phi_2=1$. Let $\{\phi_k\}_{k=3, \cdots,8}$ be a set of
generators of $H^2(T)$ such that
\[
\phi_3 \cup \phi_4=\phi_4 \cup \phi_3=\phi_5 \cup \phi_6=\phi_6 \cup
\phi_5=\phi_7 \cup \phi_8= \phi_8 \cup \phi_7=\phi_2
\]
and all other products between $\phi_k$'s, $k=3,\cdots,8$, are zero.
Let $\phi_{j+8}$ be a generator of $H^0(\{q_j\})$ for each
$j=1,\cdots,16$ such that the products in the twisted sector is
\[
\phi_k \cdot \phi_{k'} = \delta_{k',k} \phi_2
\]
for all $k=9,\cdots ,24$. The $\Q$-degree of  $\phi_k$ is $2$ for
all $k=9, \cdots ,24$, since the age of $\alpha$ on each fix points
are $1$. Let $c$ be the bijection from $\{1,\cdots,24\}$ to
$\{1,\cdots,24\}$, denoted by $k \mapsto k^c$, such that $\phi_k
\cup \phi_{k^c}=\phi_2$ for all $k\in \{1,\cdots,24\}$. The only
non-zero structure constants are $m_{1k}^k=1$ and $m_{k k^c}^2 =1$
for all $k\in \{1,\cdots,24\}$. The Euler class is $\e=2\cdot 24
\phi_2$.

Now we want to compute the multiplication on $\Z_2^I$-coinvariants
of the stringy cohomology of  the wreath product orbifold
$[T^I/\Z_2^I \rtimes \Sigma_I]$, but instead, we compute the
multiplication on the Lehn-Sorger side because of Theorem
\ref{ringiso}.

Without loss of generality, we can suppose that
$\langle\sigma,\tau\rangle$ acts transitively on $I$ and let
${gd}:={gd}(\sigma,\tau)_I$. Let $\g \in \Z_2^{o(\sigma)}$, and $\h
\in \Z_2^{o(\tau)}$ and let
\[
\bx\otimes\by :=\left(\bgot_{a\in o(\sigma)}\phi_{i_a} \right)
\otimes \left(\bgot_{b \in o(\tau)} \phi_{i_b} \right)   \in \left(
\bgot_{a \in o(\sigma)} H^{*}(T^{\g_a})^{\Z_2} \right) \otimes
\left( \bgot_{b \in o(\tau)} H^{*}(T^{\h_b})^{\Z_2}\right).
\]
The product $\bx\sigma \cdot \by\tau$ in the Lehn-Sorger algebra is
\[
\bx\sigma\cdot \by\tau = \left(\bgot_a\phi_{i_a}\right) \sigma \cdot
\left(\bgot_b \phi_{i_b}\right) \tau =
\bm_{*}\left(\prod_a\phi_{i_a} \prod_b \phi_{i_b} \cdot \e^{{gd}}
\right)\sigma\tau,
\]
where $\bm_{* }$ is the comultiplication defined in Definition
\ref{comultiplication}.

Let us observe that
\begin{eqnarray*}
\bm_{*}(\phi_1)&=&2^{|o(\sigma\tau)|-1}
\sum_{l=1}^{24}\left(\sum_{\{i_d\}=\{2,\cdots,2,l,l^c\}}
\bigotimes_d \phi_{i_d} \right), \\
\bm_{*}(\phi_2)&=&2^{|o(\sigma\tau)|-1}\bigotimes_d \phi_2, \\
\bm_{*}(\phi_k)&=&2^{|o(\sigma\tau)|-1}\sum_{\{i_d\}=\{2,\cdots,2,k\}}
\bigotimes_d \phi_{i_d},
\end{eqnarray*}
where $d \in o(\sigma\tau)$ and $k\not=1,2$. Therefore we can write
the multiplication on $H^{*}_{orb}([X/G])\{\Sigma_I\}$ as follows.
\begin{prop}
If ${gd}=0$,
\[
\bx\sigma\cdot\by\tau=
\left\{
 \begin{array}{rl}
\bm_{*}(\phi_2)\sigma\tau & \mbox{if}\ \{i_a\} \sqcup \{i_b\}=\{1,\cdots, 1, k,k^c\}\\
\bm_{*}(\phi_k)\sigma\tau & \mbox{if}\ \{i_a\} \sqcup \{i_b\}=\{1,\cdots, 1, k\}\\
 0 & \mbox{otherwise}.
 \end{array}
\right.
\]
If ${gd}=1$,
\[
\bx\sigma\cdot\by\tau=
\left\{
 \begin{array}{rl}
48 \bm_{*}(\phi_k)\sigma\tau & \mbox{if}\ \{i_a\} \sqcup \{i_b\}=\{1,\cdots, 1\}\\
0 & \mbox{otherwise}.
 \end{array}
\right.
\]
If ${gd}\geq 2$, $\bx\sigma\cdot\by\tau=0$.
\end{prop}

\section{{\bf Example : a point with the trivial
$G$-action}}\label{sec:ex2}

Let $X=pt$, a point, and let $G$ be an arbitrary finite group acting
trivially on $X$. The orbifold
cohomology of $[pt/G]$ is the center of the group ring $\C[G]$,
which is denoted by $\cZ\C[G]$. The stringy cohomology of the
trivial $(G^I \rtimes \Sigma_I)$-space $pt$ is the group ring
$\C[G^I \rtimes \Sigma_I]$. The theorem in this section suggests the
existence of the ring isomorphism between the $G^I$-coinvariants of
the stringy cohomology of $\GSI$-space $X^I$ and the Lehn-Sorger
algebra associated to $H^{*}([X/G])$ even when $G$ is not Abelian.

\begin{thm}
$\C[G_I\rtimes \Sigma_I]^{G^I}$ is isomorphic to
$\cZ\C[G]\{\Sigma_I\}$ as $\Sigma_I$-Frobenius algebras.
In particular, $\cZ\C[G]\{\Sigma_I\}$ satisfies the trace axiom.
\end{thm}

Before we start the proof, let us introduce the idempotent basis of
$\cZ \C [G]$. Let $\{ \chi_k \}_{k \in \cU}$ be the set of all
irreducible characters so that $|\cU|$ is the number of conjugacy
classes in $G$.
Define
\begin{equation}
\fu_k := \frac{\chi_k(1)}{|G|}\sum_{g \in G} \chi_k(g^{-1}) g.
\end{equation}
The following general orthogonality of characters is well-known,
c.f. \cite{I}: for all $k,l \in \cU$,
\begin{equation}\label{gen orth}
\sum_{g \in G}\chi_k(gh)\chi_l(g^{-1})=\delta_{kl} |G|
\frac{\chi_k(h)}{\chi_k(1)}.
\end{equation}
It follows that $\{\fu_k\}$ forms a basis and satisfies
\begin{equation*}
\fu_k \cdot \fu_l = \delta_{kl} \fu_k \ \text{ and } \
\eta(\fu_l,\fu_k)=\delta_{kl}\left(\frac{\chi_k(1)}{|G|}\right)^2.
\end{equation*}
In terms of this idempotent basis, $\cZ\C[G]\{\Sigma_I\}$ is
generated by $\left(\bgot_{a \in o(\sigma)} \fu_{k_a}\right) \sigma$
where $\sigma \in \Sigma_I$ and $k_a \in \cU$ for each $a \in
o(\sigma)$. By the canonical isomorphism in Proposition
\ref{linear-iso},
\begin{equation}\label{grouplinearIso}
\left(\bgot_{a \in o(\sigma)} \fu_{k_a}\right) \sigma \ \ \mapsto \
\ \frac {\prod_{a \in o(\sigma)} \chi_{k_a}(1)}{|G|^{o(\sigma)} }
\sum_{g \in G^I} \left( \prod_{a \in o(\sigma)}
\chi_{k_a}(\psi^{\sigma}(g)_a^{-1}) \right) g\sigma \ .
\end{equation}

\Pf Let $\sigma, \tau \in \Sigma_I$. We can assume that $\lan
\sigma,\tau \ran$ acts transitively on $I$ without loss of
generality. The Euler class of $\cZ\C[G]$ is $\e=\sum_{k \in \cU}
\left(\frac{|G|}{\chi_k(1)}\right)^2 \fu_k$. The product in the
Lehn-Sorger algebra in terms of the idempotent basis is given by
\begin{equation*}
\left(\bgot_{a \in o(\sigma)} \fu_{k_a}\right) \sigma \cdot
\left(\bgot_{b \in o(\tau)} \fu_{k_b}\right)\tau
\end{equation*}
\begin{equation}\label{LSprod}
=\left\{
\begin{array}{rl}
\left(\frac{|G|}{\chi_k(1)}\right)^{n+|o(\sigma\tau)|-|o(\sigma)|-|o(\tau)|}
\left(\bgot_{c\in o(\sigma\tau)} \fu_{k}\right)\cdot \sigma\tau, &
\mbox{if $k_a=k_b=k$ for all $a,b$}  \\
0, & \mbox{otherwise.}
\end{array}
\right.
\end{equation}
Hence, we need to show that
\begin{equation}\label{orbproduct}
\sum_{g \in G^I} \left( \prod_{a \in o(\sigma)}
\chi_{k_a}(\psi^{\sigma}(g)_a^{-1}) \right) g\sigma \cdot \sum_{h
\in G^I} \left( \prod_{b \in o(\tau)}
\chi_{k_b}(\psi^{\tau}(h)_b^{-1}) \right) h\tau
\end{equation}
is equal to
\begin{equation}\label{aftercontraction}
\begin{cases}
\left(\frac{|G|}{\chi_k(1)}\right)^{n} \sum_{w \in G^I} \left(
\prod_{c \in o(\sigma\tau)}
 \chi_k(\psi^{\sigma\tau}(w)_c^{-1}) \right) w\sigma\tau
 & \text{ if $k_a=k_b=k$ for all $a$ and $b$,}\\
0 & \text{ otherwise. }
\end{cases}
\end{equation}
We use the following well-known two formulas for all irreducible characters
$\chi_k, \chi_l$ and $f,h \in G$:
\begin{equation}\label{formula1}
\sum_{g \in G}
\chi_k(fg)\chi_l(g^{-1}h)=\delta_{kl}\frac{|G|}{\chi_k(1)}
\chi_k(fh)
\end{equation}
and
\begin{equation}\label{formula2}
\sum_{g \in G}
\chi_i(fghg^{-1})=\frac{|G|}{\chi_k(1)}\chi_k(f)\chi_k(h) .
\end{equation}
The first is just the generalized orthogonality of characters and
the second one is found in \cite{I}, Exercise (3.12).

Choose representatives $i_a$ and $i_b$ from $a$ and $b$ respectively
for all $a \in o(\sigma)$ and $b \in o(\tau)$. By letting $w:=
gh^{\sigma}$, we can write the expression (\ref{orbproduct}) as
\[
\sum_{w \in G^I} \cX_{w, \{k_a\}, \{k_b\}} \cdot w\sigma\tau
\]
where
\begin{equation}\label{96}
\cX_{w, \{k_a\}, \{k_b\}}:=\sum_{g \in G^I} \left[ \prod_a
\chi_{k_a}\left(g_{{i}_a}^{-1}w_{{i}_a}^{-1}\cdots
g^{-1}_{\sigma^{|a|-1}({i}_a)}w_{\sigma^{|a|-1}({i}_a)}^{-1}\right)
\cdot \prod_b\chi_{k_b}\left(g_{\sigma({i}_b)} \cdots
g_{\sigma\tau^{|b|-1}({i}_b)}\right) \right].
\end{equation}
Suppose that there is no $k \in \cU$ such that $k_a=k_b=k$ for all
$a$ and $b$. Since $\lan\sigma,\tau\ran$ acts on $I$ transitively, we
can choose our representatives $\{i_a\}$ and $\{i_b\}$ in such a way
that ${i}_{a'}=\sigma({i}_{b'})$ and $k_{a'} \not= k_{b'}$ for some $a' \in
o(\sigma)$ and $b'\in o(\tau)$. By eliminating the summation over the
component $g_{i_{a'}}$ of $g$ by Equation (\ref{formula1}), we conclude
that $\cX_{w, \{k_a\}, \{k_b\}}=0$ for each $w \in G^I$.

Now assume that $k=k_a=k_b$ for all $a$ and $b$. We need to show
that
\begin{equation}\label{9}
\cX_{w, \{k_a\}, \{k_b\}}=\left(\frac{|G|}{\chi_k(1)}\right)^{n}
\prod_{c \in o(\sigma\tau)}
 \chi_k\left(w_{{i}_c}^{-1}\cdots w_{(\sigma\tau)^{|c|-1}({i}_c)}^{-1}\right)
\end{equation}
for every $w \in G^I$.

First of all, pick a component of the summation variable $g$ in
Equation (\ref{96}) and eliminate it by using Equation
(\ref{formula1}). Next, eliminate the summation over another
component $g_i$ of $g$ by Equation (\ref{formula1}) if $g_i$ and
$g_i^{-1}$ appear in separate $\chi_k$'s, or by Equation
(\ref{formula2}) if they appear in the same $\chi_k$. Repeat this
process until all summations have been eliminated. After each
replacement, on the right-hand side of a component $g_{\sigma(i)}$
is either a $g_{\sigma\tau(i)}$ or $w_{\sigma\tau(i)}^{-1}$ and on
the left-hand side of $g_{\sigma(i)}$ is either a
$w^{-1}_{\tau^{-1}(i)}$ or $g_{\sigma\tau^{-1}(i)}$. Note that, on
the left side of $g_{\sigma(i)}^{-1}$ is always a $w^{-1}_i$. Hence,
after we eliminate the summations over all of the components of $g$
in the expression (\ref{96}), we obtain Equation (\ref{9}).

This proves that the canonical isomorphism defined in Proposition
\ref{linear-iso} preserves the ring structure. All other properties
of $\Sigma_I$-Frobenius algebras are clearly preserved by the
isomorphism. Thus $\C[G_I\rtimes \Sigma_I]^{G^I}$ is isomorphic to
$\cZ\C[G]\{\Sigma_I\}$ as $\Sigma_I$-Frobenius algebras. In
particular, $\cZ\C[G]\{\Sigma_I\}$ satisfies the trace axiom. \qed

\section{{\bf {Hilbert schemes and wreath products
orbifolds}}}\label{sec:crepant}

 In this section, we will relate the
wreath product orbifold associated to a $G$-space $X$ to the Hilbert
scheme of $n$-points on $Y$ when $Y$ is a crepant resolution of
$X/G$. Throughout the section, all vector spaces are over $\C$ and
we will work in the algebraic category.


\begin{defn}
Let $W$ be a normal variety over $\C$ and let $\sL$ be a rank $1$,
torsion free, coherent sheaf of $\sO_W$-module over $W$. $\sL$ is
called \emph{divisorial} \cite{M} if and only if any torsion
free coherent sheaf of $\sO_W$-module, $\sM$, such that $\sL \subset
\sM$ and $\Supp(\sM/\sL)$ has codimension $\geq 2$, coincides with $\sL$.
\end{defn}
\begin{rem}\label{div}
Let $\sL$ be divisorial. If $W^0 \subset W$ is a non-singular open
subvariety such that $W\backslash W^0$ has codimension $\geq 2$,
then $\sL|_{X^0}$ is invertible and $\sL=j_{\ast}(\sL|_{X^0})$
\cite{M}, where $j:W^0 \inc W$ denotes the canonical  inclusion. Let
$K_W$ be the canonical divisor of $W$. By Proposition (7) in
\cite{M}, the canonical sheaf $\omega_W:=\sO(K_W)$ of $W$ is
divisorial. Hence, we have $\omega_W=j_{*}\omega_{W^0}$ since
$\omega_W|_{W^0}=\omega_{W^0}$.
\end{rem}

\begin{defn}
Let $W$ and $Y$ be normal varieties. A birational morphism $\pi : Y
\rightarrow W$ is \emph{crepant} if $\omega_Y \cong
\pi^{*}\omega_W$.
\end{defn}

\begin{defn}\label{Gorenstein}
A normal variety $W$ is \emph{Gorenstein} if and only if all of the
local rings are Cohen-Macaulay and $K_W$ is Cartier.
\end{defn}

\begin{lem}\label{div2}
Let $W$ and $Y$ be Gorenstern varieties. If $\pi:Y\rightarrow W$ is
a birational morphism, then $\pi^{\ast}K_X$ is divisorial.
\end{lem}
\Pf Let $\dim W=\dim Y=n$. Since $K_W$ is Cartier, $\pi^{*}K_W$ is
also Cartier. Hence, $\pi_{\ast}\omega_W$ is torsion-free and of
rank 1. Let $\sM$ be a torsion-free sheaf such that $\pi^{*}\omega_W
\subset \sM$ and $\dim(\Supp \sM/\pi^{*}\omega_W) \leq n-2$. Let
$L:=K_Y-\pi^*K_W$. $L$ is Cartier and $\sL:=\sO(L)$ is an invertible
sheaf. It follows that $\pi^{*}\omega_W \otimes \sL \cong \omega_Y
\subset \sM \otimes \sL$ and
\[
\dim(\Supp ((\sM\otimes\sL)/\omega_Y))\leq \dim(\Supp
\sM/\pi^{*}\omega_W) \leq n-2.
\]
Since $H^{n-1}((\sM\otimes \sL)/\omega_Y)=0$, we have
$H^{n}(\sM\otimes \sL)\cong H^n(\omega_Y)\cong \C$ by Serre duality.
Hence, there exists an element in $\Hom(\sM\otimes\sL, \omega_Y)$
which gives a splitting of the short exact sequence
\[
0\rightarrow \omega_Y\rightarrow\sM\otimes \sL\rightarrow
(\sM\otimes\sL)/\omega_Y\rightarrow 0.
\]
However, since $\sM\otimes \sL$ is torsion-free, $(\sM\otimes
\sL)/\omega_Y=0$. \qed

\begin{thm}\label{thm:crepant}
Let $W$ and $Y$ be normal varieties with dimension $\geq 2$. Suppose
that $W\backslash W^0$ has codimension $\geq 2$ and that
$Y^n/\Sigma_n$ and $W^n/\Sigma_n$ are Gorenstein. If
$\pi:Y\rightarrow W$ is a crepant resolution, then the induced map
$\tilde{\pi} : Y^n/\Sigma_n \rightarrow W^n/\Sigma_n$ is crepant.
\end{thm}
\Pf The smooth locus of $W^n/\Sigma_n$ is equal to $(W^n\backslash
\Delta^{\star}_W)/\Sigma_n$ where $\Delta^{\star}_W$ is the pairwise
diagonal of $W^n$. Let $\cD_Y:=\pi^{-1}(\Delta^{\star}_W)$. Let
$\overline{\pi}:Y^n\backslash D_Y\rightarrow W^n\backslash
\Delta^{\star}_W$ be the map $\pi^{\times n}$ restricted to
$Y^n\backslash \cD_Y$. Since $\pi^{\times n}: Y^n \rightarrow W^n$
is crepant, $K_{Y^n}=(\pi^{\times n})^{*}K_{W^n}$.
 Consider the commutative diagram
\[
\begin{CD}
Y^n @<<< Y^n\backslash \cD_Y\\
@V\pi^{\times n} VV @V\overline{\pi}VV \\
W^n @<<< W^n\backslash \Delta^{\star}_W
\end{CD}
\]
where the horizontal arrows are the obvious inclusions. We have
\begin{equation}\label{01}
K_{Y^n\backslash \cD_Y}=K_{Y^n}|_{Y^n\backslash
\cD_Y}=\left((\pi^{\times n})^{*}K_{W^n}\right)|_{Y^n\backslash
\cD_Y}=\overline{\pi}^{*}(K_{W^n}|_{W^n\backslash
\Delta^{\star}_W})=\overline{\pi}^{*}K_{W^n\backslash
\Delta^{\star}_W}.
\end{equation}
Consider the following commutative diagram
\[
\begin{CD}
Y^n\backslash \cD_Y @>\bq>> (Y^n\backslash \cD_Y)/\Sigma_n \\
@V\overline{\pi} VV @V\tilde{\pi}'VV \\
W^n\backslash \Delta^{\star}_W @>\bq'>> (W^n\backslash
\Delta^{\star}_W)/\Sigma_n
\end{CD}
\]
where $\bq$ and $\bq'$ are the canonical projections. Since the
actions of $\Sigma_I$ on $Y^n\backslash \cD_Y$ and $W^n\backslash
\Delta^{\star}_W$ are free, Equation (\ref{01}) implies that
$K_{(Y^n\backslash
\cD_Y)/\Sigma_n}=\tilde{\pi}'^*K_{(W^n\backslash\Delta^{\star}_W)/\Sigma_n}.$
Hence
\[
K_{Y^n/\Sigma_n}|_{(Y^n\backslash \cD_Y)/\Sigma_n}=
\left(\tilde{\pi}^*K_{W^n/\Sigma_n}\right)|_{(Y^n\backslash
\cD_Y)/\Sigma_n}.
\]
Since both $K_{Y^n/\Sigma_n}$ and $\tilde{\pi}^*K_{W^n/\Sigma_n}$
are divisorial (Remark \ref{div}, Lemma \ref{div2}), we obtain
$K_{Y^n/\Sigma_n}=\tilde{\pi}^*K_{W^n/\Sigma_n}$. \qed

\begin{rem}\label{rem-goren}
For a non-singular variety $X$ with an action of a finite group $G$,
the variety $X/G$ is Gorenstein if and only if the age of $\alpha$
on any connected component is an integer for all $\alpha \in G$. See
Remark (3.2) in \cite{M}. If $\dim X$ is even and $X/G$ is
Gorenstein, by Corollary \ref{age}, $X^n/G^n\rtimes\Sigma_n$ is
Gorenstein. In particular, for a non-singular variety $Y$ with even
(complex) dimension, the age of the symmetric product $Y^n/\Sigma_n$
is always an integer so that $Y^n/\Sigma_n$ is Gorenstein.
\end{rem}

If $Y$ is a smooth projective surface, then the Hilbert-Chow
morphism $Y^{[n]} \rightarrow Y^n/\Sigma_n$ from the Hilbert scheme
of $n$ points on $Y$ to the symmetric product of $Y$ is a resolution
of singularities \cite{Fo}, which is also crepant \cite{Be}. Hence,
together with Theorem \ref{thm:crepant} and Remark \ref{rem-goren},
we obtain the following.
\begin{cor}
Let $X$ be a smooth projective surface with an action of a finite
group $G$. Suppose that $X/G$ is Gorenstein. If $\pi:Y\rightarrow
X/G$ is a crepant resolution, then $Y^{[n]} \rightarrow
W^n/\Sigma_n$ is a crepant resolution.
\end{cor}

Together with Theorem \ref{ringiso}, we obtain the following result.
\begin{thm}
Let $Y$ be a smooth projective surface with trivial canonical class. Let $X$
be a smooth projective surface with an action of a finite, Abelian
group $G$. Suppose that $X/G$ is Gorenstein. If $\pi:Y\rightarrow
X/G$ is a crepant resolution and the ordinary cohomology
ring $H^{*}(Y)$ is isomorphic as a Frobenius algebra
to the Chen-Ruan orbifold cohomology
ring $H^{*}_{orb}([X/G])$, then $Y^{[n]}\rightarrow X^n/\Sigma_n$ is
a hyper-K\"{a}hler resolution and $H^{*}( Y^{[n]})$ is isomorphic as
a ring to $H^{*}_{orb}([X^n/G^n\rtimes \Sigma_n])$.
\end{thm}
\Pf We have
\[
\cH(Y^n,\Sigma_n) \cong H^{*}(Y)\{\Sigma_n\} \cong
H^{*}_{orb}([X/G])\{\Sigma_n\} \cong \cH(X^n,G^n \rtimes
\Sigma_n)^{G^n}.
\]
where the first equality is due to \cite{FG} and the third is
Theorem \ref{ringiso}. Since $H^{*}( Y^{[n]})\cong H^{*}(Y)\{\Sigma_n\}^{\Sigma_n}$
\cite{L-S}, we obtain the theorem by taking
$\Sigma_n$-coinvariants everywhere in the above equality. \qed

This theorem is a special case of the following conjecture due to
Ruan \cite{R}.
\begin{conj}[Cohomological hyper-K\"{a}hler resolution conjecture]
Suppose that $\,Y \rightarrow X$ be a hyper-K\"{a}hler resolution of
the coarse moduli space $X$ of an orbifold $\cX$. The ordinary
cohomology ring $H^{*}(Y)$ of $\,Y$ is isomorphic to the Chen-Ruan
orbifold cohomology ring $H^{*}_{orb}(\cX)$ of $\cX$.
\end{conj}

\begin{rem}
The conjecture in the special case of wreath product orbifolds
has been verified when $X=\C^2$ and $G$ is a finite subgroup of $\SL_2(\C)$
in \cite{EG}.
In particular, an explicit ring isomorphism
between $H^*(Y^{[n]})$ and $H_{orb}^*([X^n/G^n\rtimes\Sigma_n])$
has been established when $X=\C^2$ and $G$ is a finite cyclic subgroup of $\SL_2(\C)$
by using Fock space methods in \cite{QW2}.
\end{rem}


\end{document}